\documentclass[11pt]{amsart}
\usepackage{amsmath,amsthm, amscd, amssymb, amsfonts, mathrsfs,pb-diagram,pb-xy}
\usepackage[all]{xy}
\usepackage[dvips, dvipsnames, usenames]{color}
\usepackage{enumitem}
\usepackage[section]{placeins}

\usepackage{mathrsfs}

\usepackage[T1]{fontenc}

\allowdisplaybreaks

\newcommand{\Rt}{\mathtt R}

\newcommand{\hopf}{\mathtt H}

\newcommand{\Dchaintwo}[4]{
\rule[-3\unitlength]{0pt}{8\unitlength}
\begin{picture}(14,5)(0,3)
\put(1,2){\ifthenelse{\equal{#1}{l}}{\circle*{2}}{\circle{2}}}
\put(2,2){\line(1,0){10}}
\put(13,2){\ifthenelse{\equal{#1}{r}}{\circle*{2}}{\circle{2}}}
\put(1,5){\makebox[0pt]{\scriptsize #2}}
\put(7,4){\makebox[0pt]{\scriptsize #3}}
\put(13,5){\makebox[0pt]{\scriptsize #4}}
\end{picture}}

\newcommand{\Dynkintwoxy}[5]{
\xymatrix@R-12pt{ \overset{#1}{\underset{#2}{\circ}}\ar  @{-}[r]^{#3} &
\overset{#4}{\underset{#5}{\circ}} } }

\newcommand{\acom}[1]{A_{#1}}

\newcommand{\zt}{\mathtt{z}}

\newcommand{\ot}{\otimes}

\usepackage{hhline}

\usepackage[active]{srcltx}

\newcommand{\filt}{\operatorname{filt-}}

\newcommand{\tp}{\operatorname{top}}

\newcommand{\GK}{\operatorname{GK-dim}}
\newcommand{\soc}{\operatorname{soc}}

\newcommand{\bsl}{\boldsymbol\lambda}

\newcommand{\bV}{\mathbf{V}}

\newcommand{\uf}{\mathfrak{u}}

\newcommand{\Ind}{\operatorname{Ind}}

\newcommand{\Irr}{\operatorname{Irr}}

\newcommand{\ord}{\operatorname{ord}}

\newcommand{\supp}{\operatorname{supp}}
\newcommand\ad{\operatorname{ad}}

\newcommand\lstr{\mathfrak L}
\def\bq{\mathfrak{q}}

\newcommand{\strat}{\mathfrak G}

\newcommand{\doble}{\mathfrak D}

\newcommand{\ghost}{\mathscr{G}}

\renewcommand{\_}[1]{\mbox{$_{\left( #1 \right)}$}}

\newcommand{\toba}{\mathscr B}

\newcommand{\gr}{\operatorname{gr}}
\newcommand{\grt}{\operatorname{gr-}}

\newcommand{\Fc}{{\mathcal F}}

\newcommand{\R}{{\mathcal R}}

\newcommand{\Zc}{{\mathcal Z}}

\newcommand{\ku}{\Bbbk}
\newcommand{\Gc}{{\mathcal G}}
\newcommand{\K}{{\mathcal K}}

\newcommand{\Z}{{\mathbb Z}}
\newcommand{\N}{{\mathbb N}}
\newcommand{\I}{{\mathbb I}}
\newcommand{\Ib}{{\mathbb I}}
\newcommand{\Jb}{{\mathbb J}}

\newcommand{\G}{{\mathbb G}}

\newcommand{\gi}{{\G_{\infty}}}

\newcommand{\M}{{\mathcal M}}

\newcommand{\C}{{\mathcal C}}
\newcommand{\D}{{\mathcal D}}

\newcommand{\Kc}{{\mathbf K}}

\newcommand{\Hc}{{\mathcal H}}
\newcommand{\Ac}{{\mathcal A}}

\newcommand{\Pc}{{\mathcal P}}

\newcommand{\J}{{\mathcal J}}
\newcommand{\Ic}{{\mathcal I}}

\newcommand{\ydhd}{{}^{H^*}_{H^*}\mathcal{YD}}

\newcommand{\ydhh}{{}^{H}_{H}\mathcal{YD}}
\newcommand{\ydhhf}{{}^{H}_{H}\mathcal{YD}_{\text{fd}}}

\newcommand{\ydac}{{}_{\Ac}^{\Ac}{\mathcal{YD}}}

\newcommand{\ydb}{{}_{B}^{B}{\mathcal{YD}}}

\newcommand{\ydav}{{}^{\Ac(V)}_{\Ac(V)}\mathcal{YD}}

\newcommand{\ydl}{{}_{L}^{L}{\mathcal{YD}}}

\newcommand{\ydG}{{}^{\ku G}_{\ku G}\mathcal{YD}}
\newcommand{\ydg}{{}^{\ku \Gamma}_{\ku \Gamma}\mathcal{YD}}
\newcommand{\ydk}{{}^{K}_{K}\mathcal{YD}}

\newcommand{\Ss}{\mathcal S}

\newcommand{\End}{\operatorname{End}}

\newcommand{\Rep}{\operatorname{Rep}}
\newcommand{\rep}{\operatorname{rep}}

\newcommand{\card}{\operatorname{card}}

\newcommand{\Hom}{\operatorname{Hom}}

\newcommand{\g}{{\mathfrak g}}

\newcommand{\lgot}{{\mathfrak l}}

\numberwithin{equation}{section}
\theoremstyle{plain}

\newtheorem{theorem}{Theorem}[section]
\newtheorem{lema}[theorem]{Lemma}
\newtheorem{lemma}[theorem]{Lemma}
\newtheorem{coro}[theorem]{Corollary}

\newtheorem{prop}[theorem]{Proposition}
\newtheorem{claim}{Claim}
\newtheorem{claima}{Claim}

\newtheorem{claimintro}{Claim}
\theoremstyle{definition}
\newtheorem{definition}[theorem]{Definition}

\theoremstyle{remark}
\newtheorem{obs}[theorem]{Remark}
\newtheorem{exa}[theorem]{Example}

\newtheorem{question}{Question}
\newtheorem{step}{Step}

\newcommand{\id}{\operatorname{id}}

\newcommand{\st}{\mathbb S_3}

\newcommand{\sco}{\mathbb S_5}

\def\pf{\begin{proof}}
\def\epf{\end{proof}}

\theoremstyle{remark}

\newcounter{tabla}\stepcounter{tabla}

\begin{document}

\renewcommand{\baselinestretch}{1.2}

\thispagestyle{empty}
\title[Nichols algebras over basic Hopf algebras]
{On Nichols algebras over basic Hopf algebras}

\author[Nicol\'as Andruskiewitsch and Iv\'an Angiono]
{Nicol\'as Andruskiewitsch and Iv\'an Angiono}

\thanks{The work of N. A. and I. A. was partially supported by CONICET, Secyt (UNC), the
MathAmSud project GR2HOPF}

\address{\noindent Facultad de Matem\'atica, Astronom\'{\i}a y F\'{\i}sica,
Universidad Nacional de C\'ordoba. CIEM -- CONICET. 
Medina Allende s/n (5000) Ciudad Universitaria, C\'ordoba,
Argentina}
\email{(andrus|angiono)@famaf.unc.edu.ar}


\begin{abstract} 
This is a contribution to the classification of finite-dimen\-sional Hopf algebras 
over an algebraically closed field $\ku$ of characteristic 0. Concretely, we show that a 
finite-dimensional Hopf algebra whose Hopf coradical is basic is a lifting of 
a  Nichols algebra of a semisimple Yetter-Drinfeld module
and we explain how to classify Nichols algebras of this kind.
We provide along the way new examples of 
Nichols algebras and Hopf algebras with finite Gelfand-Kirillov dimension.
\end{abstract}

\maketitle

\setcounter{tocdepth}{2}

\tableofcontents

\section{Introduction}\label{subsec:Introduction}

\subsection{The context}\label{subsec:intro-context}
We fix an algebraically closed field $\ku$ of characteristic 0.
Let $\hopf$ be a Hopf algebra with bijective antipode, $\hopf_{0}$ its coradical (the sum of its simple subcoalgebras)
and $\hopf_{[0]}$ its Hopf coradical (the subalgebra generated by $\hopf_{0}$, which is a Hopf subalgebra).
The problem of the classification of those $\hopf$ with finite Gelfand-Kirillov dimension can be organized 
in four different classes according to the following conditions, cf. \cite{AC}:

\begin{enumerate}[leftmargin=*,label=\rm{(\alph*)}]
\item\label{item:intro-class-1} $\hopf = \hopf_{0}$. That is, the class of cosemisimple Hopf algebras.

\smallbreak
\item\label{item:intro-class-2} $\hopf = \hopf_{[0]} \neq \hopf_{0}$, the class of (non-cosemisimple) Hopf algebras generated by the coradical.
	
\smallbreak
\item\label{item:intro-class-3} $\hopf \neq \hopf_{[0]} = \hopf_{0}$, i.~e. the coradical is a (proper) Hopf subalgebra.
	
\smallbreak
\item\label{item:intro-class-4} $\hopf \neq \hopf_{[0]}  \neq \hopf_{0}$.
\end{enumerate}

There is no general method, to our knowledge, for the classes \ref{item:intro-class-1} and \ref{item:intro-class-2},
even for the sub-problem of the classification of finite-dimensional Hopf algebras. For class \ref{item:intro-class-3}
there is a well-known method \cite{AS-cambr} that was applied under various natural hypothesis; see \cite{icm,AS-ann,AG} and references therein.

\smallbreak
Here we contribute to class \ref{item:intro-class-4} for the sub-problem of finite dimension, according to the method proposed in \cite{AC} that extends \cite{AS-cambr}.
Namely the graded Hopf algebra $\gr \hopf$ associated to the standard filtration splits as $R\# \hopf_{[0]}$, where $R = \oplus_{n\geq 0} R^n$ 
is a graded Hopf algebra in the category of Yetter-Drinfeld modules over $\hopf_{[0]}$ (called the \emph{diagram}). 
The proposal of \cite{AC} is to classify the possible $R$ and 
then to compute all liftings (Hopf algebra deformations) of $R\# \hopf_{[0]}$.
Our main result, Theorem \ref{thm:main},
answers fairly completely the first question under the assumption that $\hopf_{[0]}$ is basic with abelian group of characters $G$.
Namely, $R$ should be the Nichols algebra of a semisimple Yetter-Drinfeld module and the list of such Nichols algebras is controlled by the 
classification of the finite-dimensional Nichols algebras over $G$ which follows from \cite{H-classif RS}. 

The proof of Theorem \ref{thm:main} has three parts. 
The first one deals with Nichols algebras of semisimple Yetter-Drinfeld modules, see Theorem \ref{theorem:intro-basic}. 
In the second part, we prove that the Nichols algebra of a non-semisimple Yetter-Drinfeld module has infinite dimension,
see Theorem \ref{thm:intro-basic-2}. 

For the last part, we notice that in general it is not known whether an arbitrary $R= \oplus_{n\geq 0} R^n$ arising from the standard filtration is coradically graded, or whether the subalgebra $\Rt$ of $R$ generated by $\Zc = R^1$ is a Nichols algebra, albeit 
$\toba(\Zc)$ is a quotient of $\Rt$. 
Assuming that $R$ is finite-dimensional over $L$ basic with abelian group, 
we show that $R$ is a Nichols algebra in Subsection \ref{subsec:diagram}; 
the proof  relies strongly on the description of liftings of Nichols algebras over abelian groups \cite{AG,AG-survey}.

The smallest example of a non-semisimple Hopf algebra $L$ generated by its coradical 
was first considered by Radford; it has dimension 8 and is basic, 
being the dual of the pointed Hopf algebra $\mathbf r(-1)$,  
which is  a lifting of a quantum line. The study of the Nichols algebras
over $L = \mathbf r(-1)^*$ was undertaken in \cite{GGi}: our Corollary \ref{coro:gag}
generalizes \cite[Theorem A]{GGi}. We mention that one of the motivations of this paper was to put the results presented in \cite{GGi} in a general context.
Nichols algebras over other  basic, non-semisimple, Hopf algebras of small dimension were considered
in  \cite{HX,X1,X2,X3}. 

Actually, the class of basic Hopf algebras 
is a source of examples of Hopf algebras generated by its coradical, perhaps the only 
understandable presently (up to routine modifications like tensoring with a semisimple one or passing to a Morita equivalent one).
However, not every basic Hopf algebra has this property.

\begin{question}
	Given a basic Hopf algebra $L$, $\dim L < \infty$, determine when it is generated by its 
	coradical in terms of the deformation parameters of $B = L^*$.
\end{question}

\subsection{Nichols algebras over basic Hopf algebras}\label{subsec:intro-basic}

Let $L$ be a finite-dimensional \emph{basic} Hopf algebra; that is every simple $L$-module has dimension 1 or
equivalently $B = L^*$ is pointed. 
It was conjectured that such $B$ should be generated by group-like and skew-primitive elements \cite{AS-adv};
the conjecture is valid in all known examples, e.~g. when the group $G(B)$ is abelian \cite{A-presentation}.
Let $G = G(B) = \Hom_{\text{alg}} (L, \ku)$ and let $V \in \ydG$
be the infinitesimal braiding of $B$ \cite{AS-cambr}.
We assume that 
\begin{align}\label{item:intro-basic-hyp1}
&K :=\gr B \simeq \toba(V) \# \ku G,
\\ \label{item:intro-basic-hyp2}
&\text{	$B$ is a cocycle deformation of $\gr B$.}
\end{align}
Hypothesis \eqref{item:intro-basic-hyp1} is a rephrasing of the above Conjecture;  
in turn \eqref{item:intro-basic-hyp2} again holds when
$G$ is abelian \cite{AG} and in all known cases \cite{GM}.  

In some instances, we shall also need that our finite group $G$ satisfies: 
\begin{align}\label{item:intro-basic-hyp3} 
\begin{aligned}
&\text{\parbox{300pt}{\emph{Every finite dimensional pointed Hopf algebra
		with group $G$ 
is generated by group-like and skew-primitive elements}.}}
\end{aligned}
\end{align}

Evidently \eqref{item:intro-basic-hyp1}--that concerns only our $B$--is implied by \eqref{item:intro-basic-hyp3}.
By \eqref{item:intro-basic-hyp2}, there is 
an equivalence of braided tensor categories $\Fc: \ydl \to \ydk$ with inverse
\begin{align*}
\Gc: \ydk \to \ydl.
\end{align*}
Indeed, it is well-known that the Drinfeld doubles of a finite-dimensional Hopf algebra and its dual
are isomorphic, so that $\ydl \simeq \ydb$.
Since the tensor categories $B^*$-mod and $K^*$-mod are equivalent by \cite{S} and because of \eqref{item:intro-basic-hyp2},
$\ydb \simeq \ydk$ as braided tensor categories. Now, the simple objects in $\ydk$ 
are of the form $L(\lambda)$, $\lambda \in \Irr \ydG$, see  Proposition \ref{prop:L(lambda)} below.

\begin{theorem}\label{theorem:intro-basic}
Let $L$ be a basic Hopf algebra, $G$ and $V$ as above.
Assume that \eqref{item:intro-basic-hyp1}  and \eqref{item:intro-basic-hyp2} hold.
 Let $\Zc = \Gc(L(\lambda_1)) \oplus \dots \oplus \Gc(L(\lambda_t)) \in \ydl$ semisimple, 
where $\lambda_1, \dots, \lambda_t \in \Irr \ydG$. Then the following are equivalent:
\begin{align}
\label{eq:intro-basic1} \dim \toba (\Zc) &< \infty;
\\\label{eq:intro-basic2} 
\dim \toba (V \oplus \lambda_1 \oplus \dots \oplus \lambda_t) &< \infty.
\end{align}
\end{theorem}

Since both $\Fc$ and $\Gc$ preserve dimensions and Nichols algebras, the proof of
Theorem \ref{theorem:intro-basic} is reduced to the following Claim:

\begin{claimintro}\label{claimintro:a}
Let $Z = L(\lambda_1) \oplus \dots \oplus L(\lambda_t) \in \ydk$, 
where $\lambda_1, \dots, \lambda_t \in \Irr \ydG$. Then $\dim \toba (Z) < \infty$ if and only if 
\eqref{eq:intro-basic2} holds.
\end{claimintro}

Claim \ref{claimintro:a} follows from Proposition \ref{prop:B(L(lambda))}, valid for any finite-dimensional Hopf algebra $H$.
Observe that the Claim itself does not provide directly new finite-dimensional
Hopf algebras as, cf. the proof of Proposition \ref{prop:B(L(lambda))},
$$\toba(Z) \# K \simeq \toba (V \oplus \lambda_1 \oplus \dots \oplus \lambda_t) \# \ku G.$$
But the Hopf algebras of the form $\toba(\Zc) \# L$ are new, except for the small $L$ mentioned above.

\medbreak
Observe that the Nichols algebras $\toba (Z)$ bear a Weyl groupoid since $Z$ is semisimple by \cite{AHS,HS-adv}; these Weyl groupoids were studied in \cite{CL}.

\begin{theorem}\label{thm:intro-basic-2}
	Let $L$ be a basic Hopf algebra, $G$ and $V$ as above, such that $G$ is abelian.
If $\Zc \in \ydl$ has $\dim \toba (\Zc) < \infty$, then $\Zc$ is semisimple.
\end{theorem}

This is a drastic simplification, since most of the times $\ydl \simeq \ydk$ is wild.
When $G$ is abelian, Theorems \ref{theorem:intro-basic} and \ref{thm:intro-basic-2}
together with \cite{H-classif RS} reduce the complete
classification of finite-dimensional Nichols algebras in $\ydl$ to a computational problem.
Analogously to the proof of Theorem \ref{theorem:intro-basic}, Theorem \ref{thm:intro-basic-2} boils down to

\begin{claimintro}\label{claimintro:b} Assume that $G$ is a finite abelian group.
If $Z \in \ydk$ 
has $\dim \toba (Z) < \infty$, then $Z$ is semisimple.
\end{claimintro}

Claim \ref{claimintro:b} is proved as Theorem \ref{thm:Z-semisimple}.
We point out that
Theorems \ref{theorem:intro-basic} and \ref{thm:intro-basic-2} generalize, and were motivated by, \cite[Theorem 4.5]{GGi}.

\subsection{The main result}\label{subsec:intro-basic0}
In Subsection \ref{subsec:diagram} we prove that the diagram $R$ of a finite-dimensional Hopf algebra $\hopf$
whose Hopf coradical is basic is necessarily a Nichols algebra, under suitable hypothesis. Together with the results in \S \ref{subsec:intro-basic}, this rounds up the following statement.

\begin{theorem}\label{thm:main} Let $L$ be a basic finite-dimensional Hopf algebra such that $G = \Hom_{\text{alg}} (L, \ku)$ is an abelian group.
	Let $\hopf$ be a  Hopf algebra  with $\hopf_{[0]} \simeq L$, so that
	$\gr \hopf \simeq R \# L$.
	Then the following are equivalent:
	\begin{enumerate}[leftmargin=*,label=\rm{(\alph*)}]
		\item\label{item:intro-main-1} 
		$\hopf$ is finite-dimensional, i.~e. $R$ is finite-dimensional.
		
		\item\label{item:intro-main-2} 	$R \simeq \toba(\Zc)$, where $\Zc = \Gc(L(\lambda_1)) \oplus \dots \oplus \Gc(L(\lambda_t)) \in \ydl$ is semisimple, 
		with $\lambda_1, \dots, \lambda_t \in \Irr \ydG$, and 
		\begin{align}
\label{eq:intro-main2} 
		\dim \toba (V \oplus \lambda_1 \oplus \dots \oplus \lambda_t) &< \infty.
		\end{align}
	\end{enumerate}
\end{theorem}

The assumption $\hopf_{[0]} \simeq L$ contains implicitly the hypothesis that $L$ is generated by the coradical.

Theorem \ref{thm:main}
brings down the classification of the finite-dimen\-sional Hopf algebras 
with Hopf coradical $L$ (such that $G = \Hom_{\text{alg}} (L, \ku)$ is abelian) to two Questions that
we formulate in general.

\begin{question}\label{pbm:fGV} For $G$ and $V$ satisfying \eqref{item:intro-basic-hyp1} and \eqref{item:intro-basic-hyp2}, determine 
\begin{align*}
\mathfrak f_G(V) = \{U \in \ydG: U\neq 0, \ \dim \toba(U\oplus V) < \infty \},
\end{align*}
\end{question}

Since $\ydG$ is semisimple, we have a map from $\mathfrak f_G(V)$ to the class of semisimple objects in $\ydl$ as above.

\begin{question}\label{pbm:acuadra3}
	Classify all liftings of  $\toba(\Zc)\# L$, for any $\mathfrak f_G(V) \ni U \mapsto \Zc$ by the mentioned map.
\end{question}

For given $G$ and $V$,  the answer to Question \ref{pbm:fGV} follows from  \cite{H-classif RS,HV-rank>2},
up to describing the possible realizations over $G$. Here is an exhaustion result.

\begin{prop}\label{prop:fgv-empty} Let $L$, $G$, $V$ as above. 
Assume that $L$ is generated by its coradical and that  $\mathfrak f_G(V) = \emptyset$.
If $H$ is a finite-dimen\-sional Hopf algebra such that its Hopf coradical $H_{[0]} \simeq L$, then 
$H \simeq L$. \qed
\end{prop}

Indeed, if $\dim \toba(\Zc) < \infty$ for $\Zc \in \ydl - 0$, then $\dim \toba(\soc \Zc) < \infty$.

\smallbreak
In otherwords, if $\mathfrak f_G(V) = \emptyset$, then the only Hopf algebras arising from these circle of ideas are 
duals of non-trivial liftings of $\toba(V) \# \ku G$ that are generated by the coradical.

\medbreak
If $G$ is an abelian group of odd order, relatively prime to 105, then there are finitely many 
$V \in \ydG$ with $\dim \toba(V) < \infty$ \cite[Prop. 8.1]{AS-adv}, hence there are various 
$V$'s with $\mathfrak f_G(V) = \emptyset$. See also \cite[Theorem 1.3]{AS-adv} for $G \simeq \Z/p$.
Also, if $G = \st$, respectively $\sco$, and $V \in \ydG$ simple corresponding to the class of transpositions 
and the sign representation, then $\mathfrak f_G(V) = \emptyset$ by \cite[\S 4.2]{AHS}, respectively by  \cite{HV-rank>2}.

\medbreak
To implement the usual approach to  Question \ref{pbm:acuadra3}
the defining relations of the Nichols algebras $\toba\left(L(\lambda_1) \oplus \dots \oplus L(\lambda_t)\right)$
are needed; these  can be computed in principle from those of
$\toba (V \oplus \lambda_1 \oplus \dots \oplus \lambda_t)$. In Section \ref{sec:nichols-diag} we deal the case 
of diagonal type; here the existence of a PBW-basis eases up the task.
Notice that this is no enough, since the functor $\Gc$ might not be explicit.
Finally in Section \ref{sec:decomp-block} we consider decompositions involving a block coming from Nichols algebras in \cite{AAH}.

\subsection{Decompositions}\label{subsec:intro-decompositions}
The proof of Proposition \ref{prop:B(L(lambda))}
relies on a general argument allowing different variations.
Let $(W, c)$ be a braided vector space with a decomposition $W = V \oplus U$ 
such that $V$ and $U$ are braided subspaces and
\begin{align*}
c(V \otimes U) &= U \otimes V, & c(U \otimes V) &= V \otimes U.
\end{align*}
Then $\toba(W)$ is not isomorphic to 
$\toba(V) \ot\toba (U)$, unless $c_{|V \ot U}c_{|U \ot V} = \id_{U \otimes V}$ \cite{gr-jalg}. 
But there is a substitute, see \S \ref{subsec:general-coinvariant} for details. 
Namely, $\toba(W)$ splits as 
\begin{align}\label{eq:braided-bosonization-intro}
\toba(W) \simeq \K \# \toba (V),
	\end{align}
where $\K$ is an appropriate algebra of coinvariants and $\#$ stands for braided bosonization.
By \cite[Prop. 8.6]{HS-adv}, $\K$ itself is a Nichols algebra $\toba (Z_U)$ where	
\begin{align}\label{eq:braided-bosonization-intro2}
Z_U &:= \ad_c\toba (V) (U).
\end{align}
The isomorphism \eqref{eq:braided-bosonization-intro} is used in the definition of the Weyl groupoid, cf. \cite{H-inv,AHS,HS-adv}.
In \cite{AAH}, the structure of $\toba(W)$ for several families of braided vector spaces $W$ 
was determined from the knowledge of $\K$  and $\toba(V)$ via \eqref{eq:braided-bosonization-intro2}. 
In the present article, we go in the opposite direction
and get information on $\toba (Z_U)$ from the knowledge of  $\toba (W)$ and $\toba(V)$.
At least we get new examples of interesting Nichols algebras, since evidently
\begin{align} \label{eq:dimB(Z_U)}
\dim  \toba(W) &= \dim \toba (Z_U)\dim  \toba(V) ,
\\ \label{eq:GKdimB(Z_U)} \GK  \toba(W) &\leq \GK\toba (Z_U) + \GK  \toba(V),
\end{align}
with equality in \eqref{eq:GKdimB(Z_U)} if $\toba(V)$ has a convex PBW-basis, cf. \cite[Lemma 2.3.1]{AAH}
and \cite[Remark 2.3]{AAH-infinite}.
Besides,  there is some control on $Z_U$:  if $V \in \ydhh$,  $\dim H< \infty$, $\dim \toba (V) < \infty$ and
$U \in \ydhh$ is semisimple, then so is $Z_U\in {}_{\toba (V) \# H}^{\toba (V) \# H}{\mathcal{YD}}$, cf. Proposition \ref{prop:B(L(lambda))}. A closely related idea appeared in \cite[Section 4]{R} in a different 
guise and was discussed again in \cite{U}, and as mentioned above in \cite{CL}.

The paper is organized as follows. In Section \ref{sec:general} Proposition \ref{prop:B(L(lambda))} is proved. 
Section \ref{sec:indecomposable} is devoted to Nichols algebras of non-semisimple Yetter-Drinfeld modules
and the last part of the proof of the main result. 
We study Nichols algebras arising from decompositions of braided vector spaces of diagonal type in Section \ref{sec:nichols-diag}, respectively with blocks and points in Section \ref{sec:decomp-block}.
There we give partial answers to the following Question, see Theorems \ref{thm:ZU} and \ref{thm:laistrygonian}.

\begin{question}\label{question:ASreg} \cite{A-chicago}
Let $\toba(V)$ be a Nichols algebra that is a domain and has finite Gelfand-Kirillov dimension. Is $\toba(V)$ AS-regular?
\end{question}

\subsubsection*{Notations}
For us, $\N = \{1, 2, 3, \dots\}$, $\N_0 = \N \cup \{0\}$.
If $k < \theta \in \N_0$, then we denote $\I_{k, \theta} = \{n\in \N_0: k\le n \le \theta \}$,
and   $\I_{\theta} := \{1,\dots, \theta\}$; also $\I = \I_{\theta}$ if  $\theta$ is clear from the context.
The canonical basis of $\Z^{\theta}$ is denoted by $(\alpha_i)_{i\in \I_{\theta}}$.
We set
\begin{align}\label{eq:alfaij}
\alpha_{i j} &= \sum_{k \in \I_{i,j}} \alpha_k,&  i&\leq j \in \I.
\end{align}

The group of $n$-th roots of 1 in $\ku$ is denoted $\G_n$, $\G'_n$ is the subset of primitive ones, 
while $\G_{\infty} = \bigcup_{n\ge1} \G_n$, $\G'_{\infty} = \G_{\infty} - \{1\}$.

Let $\Irr \C$ be the set of isomorphism classes of irreducible objects in an abelian category $\C$. 
The category of representations of an algebra $A$ is denoted $\Rep A$ while the subcategory of finite-dimensional ones is $\rep A$.
If $M$ is a subobject of $N$ in a category $\C$, then we write $M \le N$.

See \cite{Mo-libro} for basic results and notations on Hopf algebras. 
We denote by $\Delta$ the comultiplication of a coalgebra and by $\delta$ the coaction of a comodule.
The antipode of a Hopf algebra is denoted by $\Ss$. 

Let $H$ be a Hopf algebra. As usual, $G(H)$ denotes the group of group-like elements in $H$. If $V$ is a (left) 
$H$-comodule, then $V_{g} := \{v\in V: \delta(v) = g \otimes v\}$, $g \in G(H)$.
The space of skew-primitive elements of a (braided or usual) 
Hopf algebra $H$ is denoted by $\Pc_{g,h} (H)$; that of primitive ones simply by $\Pc(H)$.

\section{Nichols algebras from decomposable braided vector spaces}\label{sec:general}

\subsection{Preliminaries}\label{subsec:general-preliminaries}
We recall the basic definitions and tools to be used along the paper.
A pair $(V, c)$ is a \emph{braided vector space} if $V$ is a vector space and $c \in GL(V \otimes V)$ satisfies
$$(c\otimes \id)(\id\otimes c)(c\otimes \id) = (\id\otimes c)(c\otimes \id)(\id\otimes c).$$ 

We assume that all Hopf algebras considered here have bijective antipode. Throughout this paper, $H$ is a Hopf algebra. 
Let $\ydhh$ be the category of Yetter-Drinfeld modules over $H$ and $\ydhhf$ the subcategory of finite-dimensional ones.
Every $V \in \ydhh$ is a (rigid) braided vector space; recall that the braiding and its inverse are
\begin{align*}
c(x \ot y) &= x\_{-1} \cdot y \ot x\_{0},& c^{-1}(x \ot y) &=  y\_{0} \ot \Ss^{-1}(y\_{-1}) \cdot x ,& x, y &\in V.
\end{align*}
Conversely every  (rigid) braided vector space
can be realized in $\ydhh$ for a suitable $H$, but by no means in a unique way.

Let $A$ be a Hopf algebra  and let  $\varkappa : H \otimes A \to \ku$ be a {\em
	skew-pairing} \cite[Definition 1.3]{DT}, i.e. a linear map such that,
for all $a,a' \in A$ and $h,h' \in H$,
\begin{align}\label{form1}
\varkappa(h,aa') &= \varkappa(h\_{2},a) \varkappa(h\_{1},a'), & \varkappa(h,1) &= \varepsilon(h),\\
\varkappa(hh',a) &= \varkappa(h, a\_{1}) \varkappa(h',a\_{2}),&
\varkappa(1,a) &= \varepsilon(a).
\end{align}
Let  $\sigma : (H \otimes A) \otimes (H \otimes A) \to k$
be the 2-cocycle associated to $\varkappa$, i.e.
\begin{align}\label{tausigma}
\sigma(h \otimes a, h' \otimes a') &= \varepsilon(h) \varkappa(h',a)\varepsilon(a'), &
a,a' \in A, \ h,h' &\in H;
\end{align}
and let $D = (H \otimes A)_{\sigma}$ be
the 2-cocycle twist of  $H \otimes A$.
Namely, $(H \otimes A)_{\sigma}$  is the tensor product coalgebra $H \otimes A$ 
with  multiplication defined by
\begin{align}\label{Eq2CocyleMult}
(h{\otimes}a)(h'{\otimes}a') &= \varkappa(h'\_{1}, a\_{1}) \varkappa^{-1}(h'_{(3)}, a_{(3)}) hh'\_{2}{\otimes}a\_{2}a',
\\
\text{where } \varkappa^{-1}(h,a) &= \varkappa(\Ss(h),a)= \varkappa(h,\Ss^{-1}(a)), 
\notag
\end{align}
$a,a' \in A$, $h,h' \in H$. The following result essentially goes back to \cite{Dr}; 
we include a proof for completeness.

\begin{lema}\label{lema:yd-rep}
	There is a tensor functor from $\ydhh$ to $\Rep D$ given by
	\begin{align*}
	\ydhh \ni M & \longmapsto M \in \Rep D,& (h \ot a) \cdot m \overset{\star}{=}  \varkappa(m\_{-1}, a) h\cdot m\_{0}, 
	\end{align*}
	$a\in A$, $h\in H$, $m \in M$.
\end{lema}

\pf Let $a,a' \in A$, $h,h' \in H$; set $h = h{\otimes}a$, $h'= h'{\otimes}a'$. Let $m \in M$. Then
\begin{align*}
&h\cdot \left(h' \cdot  m\right) 
= \varkappa(m\_{-2}, a') \varkappa(h'\_{1}m\_{-1} \Ss(h'\_{3}), a) (hh'\_{2})\cdot m\_{0}; \\
&(hh') \cdot  m  =
\varkappa(h'\_{1}, a\_{1}) \varkappa^{-1}(h'_{(3)}, a_{(3)}) \varkappa(m\_{-1}, a\_{2}a') (hh'\_{2}) \cdot  m\_{0} \\
&= \varkappa(h'\_{1}, a\_{1}) \varkappa(\Ss(h'_{(3)}), a_{(3)}) \varkappa(m\_{-2},a') \varkappa(m\_{-1}, a\_{2}) (hh'\_{2}) \cdot  m\_{0}
\end{align*}
showing that the action $\star$ is associative. By a similar computation, it is compatible with the tensor product.
\epf

Assume that $\dim H < \infty$ and take $A = (H^{\text{cop}})^*$, $\varkappa$ the standard pairing.
Then $D(H) = (H \otimes A)_{\sigma}$ is the Drinfeld double of $H$ and the functor in 
Lemma \ref{lema:yd-rep} gives an equivalence of braided tensor categories
$\ydhhf \simeq \rep D(H)$. See \cite{Majid} or \cite[10.6.16]{Mo-libro}.

\medbreak
Let $V \in \ydhh$. The left dual of $V$ is ${}^*V = \Hom(V, \ku) \in \ydhh$ with the action and coaction determined by 
\begin{align*}
\langle h\cdot f, x \rangle  &= \langle f, \Ss^{-1}(h)\cdot x \rangle, &
f\_{-1}\langle f\_{0}, x \rangle  &= \Ss (x\_{-1})\langle f, x\_{0} \rangle,
\end{align*}
$h \in H$, $x\in V$, $f \in {}^*V$. We denote by $\ad$, respectively $\ad_c$, the adjoint action of a Hopf algebra,
respectively of a Hopf algebra in $\ydhh$.
If $R$ is Hopf algebra in $\ydhh$, then the braided commutator of $x, y \in R$ is
$$[x, y]_c = xy - \text{multiplication } \circ c (x\ot y).$$

\subsection{Nichols algebras}\label{subsec:general-nichols}
Nichols algebras bring decisive information for the classification of Hopf algebras (with finite dimension or growth), see \cite{AS-cambr},
but deserve to be considered as a subject on its own, by their intricate combinatorics and potential relationship 
with other areas of Algebra. See  the survey \cite{A-leyva}. At this stage, we are interested in the following questions:

\begin{itemize}
	\item Classify all Nichols algebras with finite Gelfand-Kirillov dimension, particularly finite dimension.
	
	\item For them, provide a minimal set of defining relations.
\end{itemize}
There is no hope of a unified approach to these questions, rather one needs first to 
delimitate classes of braided vector spaces that might be approached uniformly. 
The class we understand better at this moment is that of braided vector spaces of diagonal type, see below. 
The class of braided vector spaces 
over finite groups was treated in many papers, 
with several substantial answers and lots of intriguing questions, see \cite{A-leyva}. 
The class of braided vector spaces over abelian groups but not of diagonal type was considered in \cite{AAH}.
We refer to \cite{AS-cambr,AA-diag-survey,A-leyva} for introductions to Nichols algebras,
Hopf algebras in braided categories of Yetter-Drinfeld modules, Radford-Majid theory of bosonization, etc. 
We proceed to recall the definition of Nichols algebra.

\medbreak
Let $V \in \ydhh$.
The tensor algebra $T(V)$ is a Hopf algebra in $\ydhh$ (beware, not with the usual comultiplication). 
Given $f \in {}^*V$, the skew-derivation $\partial_{f}\in \End  T(V)$  is defined by 
\begin{align}\label{eq:partial-0}
\partial_f(1) &=0,& &\partial_f(v) = f(v), \quad v\in V,\\
\label{eq:partial-right}
\partial_f(xy) &= x\partial_f(y) + \sum_j  \partial_{f_j}(x) y_j, & &\text{where } c^{-1} (y \ot f) = \sum_j f_j \ot y_j.
\end{align}
Here are some basic properties of the skew-derivations:

\begin{itemize}[leftmargin=*]\renewcommand{\labelitemi}{$\circ$}

\item The comultiplication $\Delta: T(V) \to T(V) \otimes T(V)$ is graded with respect to the standard grading of $T(V)$; 
let  $\Delta^{n-i, i}: T^n(V) \to T^{n-i}(V) \otimes T^{i}(V)$ be
its homogeneous component for $n \in\N_0$ and $i \in \I_{0,n}$.
Then an alternative definition of $\partial_{f}$ is 
\begin{align}\label{eq:partial-alt}
\partial_f &= (\id \ot f) \Delta^{n-1, 1}: T^n(V) \to T^{n-1}(V), & n \in &\N.
\end{align}

\smallbreak
\item Let $\Ic = \oplus_{n \ge 2} \Ic^n$ be a homogeneous Hopf ideal of $T(V)$ and $\R = T(V) / \Ic$. Then 
\eqref{eq:partial-alt} defines a skew-derivation $\partial_{f}\in \End  \R$, for any $f \in {}^*V$. 

\smallbreak
\item Let $f \in {}^*V$,  $h \in H$ and $x \in T(V)$. It is not difficult to see that
\begin{align}\label{eq:partial-act}
\partial_f (h \cdot x)&= h\_{1} \cdot \partial_{\Ss(h\_{2}) \cdot f} (x).
\end{align}

\smallbreak
\item 
Suppose that there exist a basis $(x_i)_{i\in \I}$ of $V$  and 
a family $(g_i)_{i\in \I}$ in $G(H)$ such that $\delta(x_i) = g_i \ot x_i$, for $i \in \I$.
Let $(f_i)_{i\in \I}$ be the dual basis to the previous one and $\partial_i = \partial_{f_i}$, $i\in \I$. 
Then to require \eqref{eq:partial-right} for all $f$ is equivalent to require that for all $i\in \I$
\begin{align}\label{eq:partial-right-diag}
\partial_i(xy) &= x\partial_i(y) + \partial_i(x) \,g_i\cdot y,&
x,y&\in T(V).
\end{align}
\end{itemize}

\begin{definition}
Let $\J(V) = \oplus_{n \ge 2}  \J^n(V)$, where for $ n \geq 2$,
\begin{align}\label{eq:nichols-def-der}
\J^n(V) &= \{x \in T^n(V):\, \forall f_1, \dots, f_n \in {}^*V,  \partial_{f_1} \dots \partial_{f_n} (x) =0\}.
\end{align}
Then $\toba(V) := T(V)/ \J(V)$ is a Hopf algebra in $\ydhh$, called the Nichols algebra of $V$; see e.g. \cite{AHS} for more details on this approach. 
\end{definition}

Notice that as an algebra (and coalgebra) $\toba(V)$ does depend only on the braiding
$c$, by definition of the skew-derivations $\partial_f$. 
As observed, the skew-derivations can be extended to $\toba(V)$ and by definition, we have 
\begin{align}\label{eq:deriv-nichols}
\bigcap_{f\in {}^*V} \ker \partial_f &=\ku 1 \text{ in } \toba(V).
\end{align}

\subsection{The coinvariant Nichols algebra}\label{subsec:general-coinvariant}
We explain how Nichols algebras arise from decomposition, a crucial point in this paper.
Let $H$ be a Hopf algebra. Let $V, U\in \ydhh$ and 
\begin{align*}
W = V \oplus U.
\end{align*}  
Thus we have a decomposition of $W$ as in \S \ref{subsec:intro-decompositions}. 
Conversely, any decomposition as in \S \ref{subsec:intro-decompositions}
can be realized over a suitable $H$ provided that the braiding of $W$ is rigid.

For simplicity, we set $\Ac(V) = \toba (V)\#  H$, etc.
The  natural morphisms  
(of Hopf algebras in $\ydhh$) $\toba (W)  \to \toba (V)$ and $\toba (V)  \to \toba (W)$ 
induce--by tensoring with $\id_H$--morphisms of Hopf algebras 
\begin{align*}
\pi_{\Ac(V)}: \Ac(W) &\to \Ac(V),& & \text{and}&  \iota_{\Ac(V)}: \Ac(V) &\to \Ac(W).
\end{align*}
Since $\pi_{\Ac(V)}\iota_{\Ac(V)} = \id_{\Ac(V)}$, the Radford-Majid theory applies. Namely,
$$\K = \Ac(W)^{\mathrm{co}\,\pi_{\Ac(V)}}$$ is a Hopf algebra in ${}^{\Ac(V)}_{\Ac(V)}\mathcal{YD}$ with the adjoint action
and the coaction 
\begin{align} \label{eq:coaction-R}
\delta &=(\pi _{\Ac(V)}\otimes \id)\Delta _{\Ac(W)}.
\end{align}

\begin{obs} Let $\pi_H: \Ac(V) \to H$ and $\iota_H: H \to \Ac(V)$ be the natural projection and inclusion respectively.
They give rise to restriction and corestriction functors ${}_{\Ac(V)} \M \to {}_{H} \M$
and ${}^{\Ac(V)} \M \to {}^{H} \M$; it is easy to see that they glue together to a tensor functor	
${}^{\Ac(V)}_{\Ac(V)}\mathcal{YD} \to \ydhh$, that however does not preserve braidings.
\end{obs}

Furthermore, $\Ac(W)$ is the biproduct or bosonization of $\K$, that is
\begin{align*}
\Ac(W) & \simeq \K \# \Ac(V).
\end{align*}
In fact,  $\toba(W)$ is the braided bosonization $\K \# \toba (V)$, see \cite{AHS,HS-adv} for details.
The next result generalizes \cite[Proposition 22]{R}.

\begin{prop}\label{prop:braided-bosonization} \cite[Proposition 8.6]{HS-adv}.
	$\K\simeq \toba(Z_U)$, where
	\begin{align}\label{eq:braided-bosonization}
Z_U &:= \ad_c\toba (V) (U)  \in {}_{\Ac(V)}^{\Ac(V)}{\mathcal{YD}}.
\qed
\end{align}
\end{prop}
Of course, interchanging $U$ and $V$, we obtain another Nichols algebra $\toba (Z_V)$, where
$Z_V:= \ad_c\toba (U) (V) \in {}_{\Ac(U)}^{\Ac(U)}{\mathcal{YD}}$.

\begin{obs}
If $c_{V\ot U}c_{U\ot V} = \id_{U\ot V}$, then $Z_U = U$ and $\K\simeq \toba(U)$.
\end{obs}

\bigbreak
Notice that $W$ and consequently $\toba(W)$ and $\K$ are $\Z^{2}$-graded  by 
\begin{align*}
\deg V &= \alpha_{1},& \deg U &= \alpha_2,
\end{align*}
 cf. \cite{AHS}. Hence
every pair of integers $(d_1, d_2)$ defines a $\Z$-grading of $\toba(W)$ 
via the group homomorphism $\Z^{2} \to \Z$, $\alpha_i \mapsto d_i$. 
The usual $\Z$-grading of $\toba(W)$ comes from the pair $(1,1)$.
Also, the grading of $\K$ as a Nichols algebra  arises from the given by $(0,1)$, i.e.
$\deg U = 1$, $\deg V = 0$; thus $\deg Z_U = 1$.

Consider the grading from the pair $(1,0)$. Then
$Z_U = \oplus_{j\geq 0} Z^j_U$ is $\Z$-graded, with
\begin{align*}
Z^0_U  &= U ,& Z^j_U  &= \ad_c\toba^j (V) (U) \leq \toba^{j + 1} (W) \text{ in } \ydhh.
\end{align*}
\begin{lema}\label{lem:graded-simple} Let $N \leq Z_U$ in $ {}_{\Ac(V)}^{\Ac(V)}{\mathcal{YD}}$.
If $N = \oplus_{j\geq 1}(N \cap  Z^j_U)$, then $N = 0$.
\end{lema}

\pf Let $\Ic$ be the 2-sided ideal of $\toba(W)$ generated by $N$. We claim that $\Ic$ is a Hopf ideal.
Indeed, let $\varpi_{\toba(W)}: \Ac(W) \to \toba(W)$ be the usual projection. Write $\delta(n) = n\_{-1} \ot n\_{0} \in \Ac(V) \ot N$ and compute
\begin{align*}
\Delta_{\Ac(W)} (n) &= n \ot 1 + n\_{-1} \ot n\_{0}  \implies\\
\Delta_{\toba (W)} (n) &= n \ot 1 + \varpi_{\toba(W)} (n\_{-1}) \ot n\_{0} \in N\ot \toba(W) + \toba(W) \ot N, 
\end{align*}
since $N$ is a subcomodule of $Z_U$. 
Now, if $M$ is a subcomodule of $Z_U$, then
\begin{align*}
c(M \ot \toba(W)) &= \toba(W) \ot M,& c(\toba(W) \ot M) &= M \ot \toba(W),
\end{align*}
braiding in $\ydhh$, since the action and the coaction of $H$ on $M$ are the restriction of the action of $\Ac(V)$,
respectively the corestriction of the coaction of $\Ac(V)$. The claim follows by a recursive argument.
Since $N$ is positively graded by hypothesis,
$\Ic \subseteq \oplus_{n \ge 2}\toba^n(W)$.  Thus $\Ic$, and a fortiori $N$, is 0.
\epf

\smallbreak
We next characterize $\toba(Z_U)$ using the skew-derivations $\partial_f$. Clearly, we have a decomposition ${}^*W \simeq {}^*V \oplus {}^*U$ in $\ydhh$,
thus we have skew-derivations $\partial_{f}$ of $\toba(W)$ for all $f\in {}^*V$, extending those of $\toba(V)$.

\begin{prop}\label{prop:B(ZU)-derivations}
	$\toba(Z_U)= \bigcap_{f\in {}^*V} \ker \partial_f$.
\end{prop}
\pf We first claim that $\Kc := \bigcap_{f\in {}^*V} \ker \partial_f$ is a subalgebra of $\toba(W)$; this follows from 
\eqref{eq:partial-right} since $c^{-1} (\toba(W) \ot {}^*V) = {}^*V \ot \toba(W)$.
Also, $\Kc$ is stable under the action of $H$ by \eqref{eq:partial-act}. We next claim that $\ad_c(V)(\Kc) \subseteq \Kc$.
For, let $f\in {}^*V$, $x \in V$ and $y \in \Kc$. Then
\begin{align*}
\partial_{f}(\ad_c(x) y) &= \partial_{f}\left(xy - (x\_{-1} \cdot y) x\_{0} \right)
 = x\partial_{f}(y) + \langle f\_{0}, x\rangle  \Ss^{-1} (f\_{1}) \cdot y
\\ & - (x\_{-1} \cdot y) \langle  f, x\_{0}\rangle - \partial_{f\_{0}}(x\_{-1} \cdot y) \Ss^{-1}(f\_{1}) \cdot x\_{0} =0.
\end{align*}

Now $U\subset \Kc$, hence 
$Z_U = \ad_c\toba (V) (U) \subset \Kc$, thus
$\toba(Z_U) \subseteq  \Kc$, since $\toba(Z_U)$ is generated by $Z_U$ as algebra. 

Conversely, let $x \in\Kc$. Fix a basis $(k_s)_{s \in S}$ of $\toba(Z_U)$. 
As the multiplication induces a linear isomorphism
$\toba(W) \simeq \toba(Z_U) \# \toba (V)$, we have
\begin{align*}
x&= \sum_{s \in S} k_s a_s, &   \text{for some } a_s & \in \toba (V).
\end{align*}
Given $f\in {}^*V$, we have that $0=\partial_f(x)= \sum_{s\in S} k_s \partial_f(a_s)$ since $k_s\in\ker \partial_f$, 
thus $\partial_f(a_s)=0$, hence $a_s \in \ku$ for all $s \in S$ by \eqref{eq:deriv-nichols}. Thus  $x\in \toba(Z_U)$.
\epf

\subsection{Nichols algebras of semisimple Yetter-Drinfeld modules}\label{subsec:ydm}

Let $H$ be a Hopf algebra with bijective antipode,  $V\in \ydhh$ and
$\Ac(V) = \toba(V) \# H$ as before. 
We first recall the construction of a family of irreducible Yetter-Drinfeld modules of $\Ac(V)$ 
from \cite[\S 3.3]{AHS}, cf. also \cite{R}.
Given $\lambda \in \Irr \ydhh$, 
take $W = V\oplus \lambda \in \ydhh$, so that 
$\toba(W)$ has a decomposition $\toba(W)=\K\#\toba(V)$. Then $L(\lambda) := \ad_c \toba(V)(\lambda) \subset \K$
belongs to $\ydav$.

\begin{prop}\label{prop:L(lambda)-YD} \cite[Lemma 3.3 \& Proposition 3.5]{AHS}
The module $L(\lambda)$ is simple and the map $L:\Irr \ydhh \to \Irr \ydav$, $\lambda \mapsto L(\lambda)$,
is injective. 
Also $L(\lambda)= \oplus_{n\in\N} L(\lambda)^n$ is a graded Yetter-Drinfeld submodule of $\Pc(\K)$ over $\Ac(V)$ (and over $H$)  with $L(\lambda)^1=\lambda$. \qed
\end{prop}

Under some conditions on $H$ and $V$  the family of simple Yetter-Drinfeld modules $L(\lambda)$ exhaust those of
$\ydav$ or of a closely related category. For simplicity of the exposition, 
and because of our main goals in this paper,
we shall consider the following  assumptions:
\begin{flalign}
\label{eq:hypothesis-1} &H \text{ is a finite-dimensional Hopf algebra.}
\\
\label{eq:hypothesis-2} &V \in \ydhh \text{ has }\dim \toba(V) < \infty.
\end{flalign}
If both assumptions hold, then $\Ac(V)$ has finite dimension and its dual is isomorphic to $\toba(\overline{V}) \# H^*$, where 
$\overline{V} = \Hom(V, \ku) \in \ydhd$ in such a way that for $v\in V$, $x \in \overline{V}$,
$h\in H$, $a \in H^*$, 
\begin{align}\label{eq:Vdual-ydhd}
(x| h\cdot v) &= (x\_{-1} | h) (x\_{0} | v),& 
(a \cdot x| v) &= (a | \Ss(v\_{-1})) (x| v\_{0}).
\end{align}
Let $\doble = D(\toba(V) \# H)$ be the Drinfeld double of $\toba(V) \# H$.
The multiplication induces a triangular decomposition: 
\begin{align*}
\doble &\simeq \toba (V) \ot D(H) \ot \toba (\overline{V}).
\end{align*}
Now the assignment $\deg V = 1 = -\deg \overline{V}$, $\deg D(H) =0$, induces a grading on $\doble$; 
notice that this is opposite to the grading in \cite{PV}. Accordingly,
\begin{align*}
\doble^{\le 0} := \toba (\overline{V}) \# D(H) \simeq   D(H) \ot \toba (\overline{V}) \hookrightarrow \doble \hookleftarrow \doble^{\ge 0} := \toba (V) \# D(H).
\end{align*}

Every $D(H)$-module becomes a $\doble^{\le 0}$-module with the trivial action of $\toba (\overline{V})$.
Let $\lambda \in \Irr \Rep D(H) \simeq \Irr \ydhh$. We consider  the Verma module
\begin{align*}
M(\lambda) = \Ind_{\doble^{\le 0} }^\doble \lambda = \doble \otimes _{\doble^{\le 0} } \lambda \simeq \toba (V) \ot \lambda.
\end{align*}
Notice that $M(\lambda) = \oplus_{n\ge 0} M^n(\lambda)$ is a graded $\doble$-module, where $M^n(\lambda) = \toba^n (V) \ot \lambda$ up to the isomorphism above.
Thus $M^0(\lambda) \simeq \lambda$ in $\ydhh$.
Set $M_+(\lambda) = \oplus_{n > 0} M^n(\lambda)$.

The following result is known for $H = \ku G$ a group algebra; it
was proved in \cite{HY-Shap-det} when  $G$ abelian and in \cite{PV} for any finite $G$,
although the main ideas of the proof appeared much earlier \cite{Lu,RS}.
Here we only assume \eqref{eq:hypothesis-1} and \eqref{eq:hypothesis-2}.
We consider the simple modules $L(\lambda)$ of Proposition \ref{prop:L(lambda)-YD} with the grading shifted in $-1$; that is, $L(\lambda)= \oplus_{n\ge 0} L(\lambda)^n$, where $L(\lambda)^0=\lambda$.

\begin{prop}\label{prop:L(lambda)}
The Verma module $M(\lambda)$ has a unique simple quotient, which is graded and 
isomorphic to $L(\lambda)$.
The family $L(\lambda)$,  $\lambda \in \Irr \ydhh$, is a parame\-trization of $\Irr \Rep \doble$.	
\end{prop}

\pf We start by the uniqueness; cf. the proofs of \cite[Theorems 1 \& 3]{PV}.

\begin{claim}\label{claim:1}
	Let $x \in M^0(\lambda) \simeq \lambda$, $x \neq 0$ and $y \in M_+(\lambda)$. Then 
	$M(\lambda) = \doble^{\ge 0} (x + y)$.
\end{claim}

Indeed, let $N =  \doble^{\ge 0} (x + y)$. If $y =0$, then $D(H) \cdot x = M^0(\lambda) \hookrightarrow N$, since $\lambda$ is simple,
and $\doble^{\ge 0} \cdot x = \toba(V) \cdot M^0(\lambda) = N$.
Assume that $y \neq 0$ and write $y = y_1 + y_2$, where $0 \neq y_1 \in  M^{h}(\lambda)$ and $y_2 \in \oplus_{j > h} M^{j}(\lambda)$. 
By the preceding, there exists $z \in \toba^{h}(V) \# D(H)$ such that $z\cdot x = y_1$. Then
\begin{align*}
x + y - z\cdot (x + y) &= x + y_2 - z\cdot y \in N, & &\text{ and }&  y_2 - z\cdot y &\in \oplus_{j > h} M^{j}(\lambda).
\end{align*}
Arguing recursively, we may assume that $ h = \tp$ but then $y_2 - z\cdot y = 0$. 

\begin{claim} $M(\lambda)$ has a unique maximal submodule
	$N(\lambda)$, which is graded.
\end{claim}

Indeed, let $N(\lambda)$ be the sum of all graded submodules contained in $M_+(\lambda)$.
Let $N \le M(\lambda)$ and let $\widetilde{N}$ be the span of all homogeneous components of all elements in $N$;
clearly $\widetilde{N}$ is also a submodule of $M(\lambda)$. Now $N \neq M(\lambda)$ implies $\widetilde{N} \subseteq M_+(\lambda)$
by Claim \ref{claim:1}. Thus any proper submodule of $M(\lambda)$ is contained in $N(\lambda)$, as expected.

Thus we have a unique simple quotient  $M(\lambda)/N(\lambda)$, which is graded.
Let us see that it is isomorphic to $L(\lambda)$.
By Lemma \ref{lema:yd-rep} we may consider $L(\lambda)$ as a (graded) $\doble$-module. 
Hence $\overline{V}$ acts trivially on $L(\lambda)^0=\lambda$, by looking at the degree.
Hence there exists a map $\pi:M(\lambda)\to L(\lambda)$ of $\doble$-modules; 
as $L(\lambda)$ is generated by $\toba(V)$ as $\toba(V)\# H$-module, $\pi$ is surjective.
Since $L(\lambda)$ is simple by Proposition \ref{prop:L(lambda)-YD}, we have $L(\lambda) \simeq M(\lambda)/N(\lambda)$.

Finally, let $L$ be a finite-dimensional $\doble$-module and let $S$ be a simple $\doble^{\ge 0}$-submodule of $L$. 
Necessarily $ S \simeq \lambda$ for some $\lambda \in \Irr \ydhh$ with trivial action of $\toba(V)$.
Thus we have a non-zero morphism of $\D$-modules $M(\lambda) \to L$; if $L$ is simple, then this is an epimorphism
and thus $L \simeq L(\lambda)$.
\epf

\bigbreak
We have identified the simple modules from Proposition \ref{prop:L(lambda)-YD} as quotients of Verma modules, 
but we can also derive Proposition \ref{prop:L(lambda)-YD} from Lemma \ref{lem:graded-simple}. Say for this purpose that 
$L(\lambda)$ is the unique simple quotient of $M(\lambda)$.

\begin{prop}\label{prop:B(L(lambda))}
	Let $\lambda_1, \dots\lambda_t \in \Irr \ydhh$. Then 
	\begin{align*}
	\toba(V \oplus \lambda_1 \oplus \dots \oplus\lambda_t) \simeq \toba(L(\lambda_1) \oplus \dots \oplus L(\lambda_t))\# \toba(V).
	\end{align*}
\end{prop}

\pf In the context of Proposition \ref{prop:braided-bosonization}, take $U= \lambda_1 \oplus \dots \oplus\lambda_t$. 
We have a projection $\pi:  M(\lambda_1) \oplus \dots \oplus M(\lambda_t) \to Z_{\lambda_1 \oplus \dots \oplus\lambda_t}$, 
homogeneous with respect to
the grading considered in Lemma \ref{lem:graded-simple}, i.~e. corresponding to $(1,0)$. 
Let  $N = \pi(N(\lambda_1) \oplus \dots \oplus N(\lambda_t))$;
clearly $N \subset \oplus_{j> 0} Z^j_U$. By Lemma \ref{lem:graded-simple}, we conclude that $N=0$.
\epf

In this way, Propositions \ref{prop:L(lambda)} and \ref{prop:B(L(lambda))} 
reduce the calculation of Nichols algebras of semisimple Yetter-Drinfeld modules
over $\Ac(V)$ to the knowledge of Nichols algebras of  semisimple Yetter-Drinfeld modules over $H$.
As an illustration, we work out the example considered in \cite{GGi}.

\begin{exa}\label{exa:gag}
Let $\Gamma \simeq \Z /N$ where $1 < N \in \N$. Fix $g \in \Gamma$ a generator, 
so that $\Gamma = \{g^a: a \in \Z/N\}$; fix $q \in \G'_N$ and let $\eta \in \widehat{\Gamma}$
given by $\eta(g) = q$, so that $\widehat{\Gamma} = \{\eta^b: b \in \Z/N\}$.
Let $\ku^{b}_{a} = \ku^{\eta^b}_{g^a} \in \ydg$ be the one-dimensional object with action given by $\eta^b$
and coaction by $g^a$. If $a,b,c,d \in \Z/N$, then $\ku^{b}_{a} \oplus \ku^{d}_{c}$ is of diagonal 
type with Dynkin diagram 
$\xymatrix@R-12pt{ \overset{q^{ab} }\circ\ar  @{-}[r]^{q^{ad+bc}}  & \overset{q^{cd}}\circ }$.
\end{exa}

Assume that $N = 4$ and let $V = \ku^{2}_{1}$, so that $\toba(V) = \Lambda(V)$ has dimension 2. 
Let $\Theta = \Theta_0 \cup \Theta_1 \cup \Theta_3$, where
\begin{align*}
\Theta_0 &= \{(1,2), (3,2)\}, \, \Theta_1 = \{(1,3), (2, 1),  (3,3)\},\, \Theta_3 = \{ (1,1), (2,3), (3, 1)\}.
\end{align*}
If $(a,b) \in \Z /N \times \Z /N - \Theta$, then $\dim \toba(\ku^{b}_{a}) = \infty$.
The following claims are proved by inspection in the classification list \cite{H-classif RS}.

\begin{claima}\label{claima:uno} The following are equivalent:
\begin{enumerate}
\item $\dim \toba (V \oplus \ku^{b}_{a}) < \infty$.

\item $(a,b) \in \Theta - \{(3,1), (3,3)\}$.
\end{enumerate}
If this is the case, then:

\begin{itemize}
	\item For  $(a,b) \in \Theta_0$, we have $\dim \toba (Z_{\ku^{b}_{a}}) = 2$. Indeed the Dynkin diagram of $V \oplus \ku^{b}_{a}$ is in both cases $\xymatrix@R-12pt{ \overset{-1 }\circ    & \overset{-1}\circ }$.
	
	\item For  $(a,b) \in \{(1,3), (2, 1), (1,1), (2,3)\}$, we have $\dim \toba (Z_{\ku^{b}_{a}}) = 8$.
	 Indeed the Dynkin diagrams of $V \oplus \ku^{b}_{a}$ are respectively
\begin{align*}
&\xymatrix@R-12pt{ \overset{-1 }\circ\ar  @{-}[r]^{q}  & \overset{q^{-1}}\circ },&
&\xymatrix@R-12pt{ \overset{-1 }\circ\ar  @{-}[r]^{q}  & \overset{-1}\circ },&
&\xymatrix@R-12pt{ \overset{-1 }\circ\ar  @{-}[r]^{q^{-1}}  & \overset{q}\circ },&
&\xymatrix@R-12pt{ \overset{-1 }\circ\ar  @{-}[r]^{q^{-1}}  & \overset{-1}\circ }.
\end{align*}
\end{itemize}
These are of type $\mathbf A(1\vert 1)$, see \cite[\S 5.1.11]{AA-diag-survey}, and $\dim \toba (V \oplus \ku^{b}_{a}) = 16$.
\end{claima}

\begin{claima} Let  $(a_i,b_i) \in \Theta_0$, $i\in \I_M$, and $U =\oplus_{i\in \I_M} \ku^{b_i}_{a_i}$. Then
$Z_{U}  = U$ and $\toba(U)$ has dimension $2^M$. Indeed the Dynkin diagram of $V \oplus U$ 
consists of disconnected points.	
\end{claima}

\begin{claima}\label{claima:tres} 
Let $(a,b), (c,d) \in \Theta - \{(3,1), (3,3)\}$ with $(a,b) \notin \Theta_0$, and $U = \ku^{b}_{a} \oplus \ku^{d}_{c}$. 
Then the following are equivalent:
\begin{enumerate}
	\item $\dim \toba (V \oplus U) < \infty$.
	
	\item  One of the next possibilities occurs:
\end{enumerate}

\begin{enumerate}[leftmargin=*,label=\rm{(\alph*)}]
	\item  $U = \ku^{3}_{1} \oplus \ku^{2}_{1}$;  the Dynkin diagram of $V \oplus U$ is 
	$\xymatrix@R-12pt{ \overset{-1 }\circ \ar @{-}[r]^{q} & \overset{q^{-1}}\circ \ar  @{-}[r]^{q}  & \overset{-1}\circ }$.

\medbreak
	\item  $U = \ku^{3}_{1} \oplus \ku^{1}_{1}$;  the Dynkin diagram of $\ku^{3}_{1} \oplus V \oplus \ku^{1}_{1}$ is 
	$\xymatrix@R-12pt{ \overset{q^{-1}}\circ \ar @{-}[r]^{q} & \overset{-1 }\circ \ar  @{-}[r]^{q^{-1}}  & \overset{q}\circ }$. 

\medbreak
\item  $U = \ku^{1}_{2} \oplus \ku^{2}_{3}$;  the Dynkin diagram of $V \oplus U$ is 
	$\xymatrix@R-12pt{ \overset{-1 }\circ \ar @{-}[r]^{q} & \overset{-1}\circ \ar  @{-}[r]^{q^{-1}}  & \overset{-1}\circ }$.

\medbreak
\item  $U = \ku^{1}_{2} \oplus \ku^{3}_{2}$;  the Dynkin diagram of $\ku^{1}_{2} \oplus V \oplus \ku^{3}_{2}$ is 
$\xymatrix@R-12pt{ \overset{-1}\circ \ar @{-}[r]^{q} & \overset{-1 }\circ \ar  @{-}[r]^{q^{-1}}  & \overset{-1}\circ }$.

\medbreak
\item  $U = \ku^{1}_{1} \oplus \ku^{2}_{1}$;  the Dynkin diagram of $V \oplus U$ is 
$\xymatrix@R-12pt{ \overset{-1 }\circ \ar @{-}[r]^{q^{-1}} & \overset{q}\circ \ar  @{-}[r]^{q^{-1}}  & \overset{-1}\circ }$.

\medbreak
 \item  $U = \ku^{3}_{2} \oplus \ku^{2}_{3}$;  the Dynkin diagram of $V \oplus U$ is 
$\xymatrix@R-12pt{ \overset{-1 }\circ \ar @{-}[r]^{q^{-1}} & \overset{-1}\circ \ar  @{-}[r]^{q}  & \overset{-1}\circ }$.	
\end{enumerate}	

Now $\dim \toba (V \oplus U) = 256$ because all these Dynkin diagrams are of 
type $\mathbf A(2\vert 1)$ \cite[\S 5.1.8]{AA-diag-survey}.
Thus $\dim \toba(Z_U) = 128$ in all cases.		
\end{claima}

\begin{claima}
Let $U\in \ydg$ of dimension 3. Then either $\dim \toba(Z_U) = \infty$ or all simple submodules of $U$ have 
labels in $\Theta_0$.
\end{claima}

Indeed,  $U = \ku^{b}_{a} \oplus \ku^{d}_{c}\oplus \ku^{f}_{e}$, where 
$(a,b), (c,d), (e,f) \in \Theta - \{(3,1), (3,3)\}$. We may assume that $(a,b) \notin \Theta_0$.
Then $\ku^{b}_{a} \oplus \ku^{d}_{c}$ and $\ku^{b}_{a} \oplus \ku^{f}_{e}$ belong to the list in Claim \ref{claima:tres},
reducing drastically the possibilities. Then we proceed by inspection again.
Collecting together these Claims, we have:

\begin{prop}\label{prop:radford1}
Let $K = \toba(V) \# \ku \Gamma$, where $\Gamma$ is cyclic of order 4 and $V = \ku^{2}_1$ as above.
Let $Z = \oplus_{i\in \I_M} L(\ku^{b_i}_{a_i})\in \ydk$. Then the following are equivalent:

\begin{enumerate}
	\item $\dim \toba (Z) < \infty$.
	
	\item  One of the next possibilities occurs:
\end{enumerate}

\begin{enumerate}[leftmargin=*,label=\rm{(\alph*)}]
\item  All $L(\ku^{b_i}_{a_i})$ have dimension 1; then $\dim \toba (Z) = 2^M$.

\medbreak
\item $M =1$ and $\ku^{b}_{a}$ is as in  Claim \ref{claima:uno}. There are 4 examples, all with $\dim \toba (Z) = 8$.

\medbreak
\item $M =2$ and $\ku^{b_1}_{a_1} \oplus \ku^{b_2}_{a_2}$ is as in  Claim \ref{claima:tres}.
There are 6 examples, all with $\dim \toba (Z) = 128$.
\end{enumerate}\qed
\end{prop}

Let $K$ be as in Proposition \ref{prop:radford1}, let $L$ be the Radford Hopf algebra of dimension 8 and let $\Gc: \ydk \to \ydl$ be the equivalence of braided tensor categories as in the Introduction. 
The classification of the finite-dimensional Nichols algebras $\toba(Z)$ with $Z \in \ydl$ was addressed in \cite{GGi}. 
If $\dim \toba(Z) < \infty$, then $Z$ should be semisimple by \cite[Theorem 4.5]{GGi} and
the classification is achieved assuming that $Z$ is simple \cite[Theorem A]{GGi}. 
The classification might be concluded as an application of the previous result. 

\begin{coro}\label{coro:gag}
Let $Z \in \ydl$. The following are equivalent:
	
	\begin{enumerate}
		\item $\dim \toba (Z) < \infty$.
		
		\item $Z$ is semisimple and one of the next possibilities occurs:
	\end{enumerate}
	
	\begin{enumerate}[leftmargin=*,label=\rm{(\alph*)}]
		\item  All simple submodules of $Z$ have dimension 1; then $\dim \toba (Z) = 2^{\dim Z}$.
		
		\medbreak
		\item $Z$ is one of the simple objects $\Gc(L(\ku^{b}_{a}))$, where  $\ku^{b}_{a}$ is as in  Claim \ref{claima:uno}.
		There are 4 examples, all with $\dim \toba (Z) = 8$.
		
		\medbreak
		\item $Z$ is one of $\Gc(L(\ku^{b_1}_{a_1})) \oplus \Gc(L(\ku^{b_2}_{a_2}))$
		where $\ku^{b_1}_{a_1} \oplus \ku^{b_2}_{a_2}$ is as in  Claim \ref{claima:tres}.
		There are 6 examples, all with $\dim \toba (Z) = 128$.
	\end{enumerate}
\end{coro}

To be precise one needs to identify the simple objects $\Gc(L(\ku^{b}_{a}))$.

\section{Semisimplicity and the diagram} \label{sec:indecomposable}

Here we prove the second and third parts of the proof of our main Theorem \ref{thm:main}.

\subsection{Nichols algebras of graded or filtered Yetter-Drinfeld modules}\label{subsubsec:graded-yd}
We start by some generalities mostly well-known. 

\subsubsection{Graded} Let $\Ac = \oplus_{i\in \Z} \Ac^i$ be a graded Hopf algebra, i.e. $\Ac^i\Ac^j \subseteq \Ac^{i+j}$ and 
$\Delta(\Ac^i) \subseteq \sum_{h + k = i} \Ac^h \ot \Ac^k$ for all $i, j \in \Z$ (then the antipode preserves the grading).

A {graded} Yetter-Drinfeld module over $\Ac$ is $M \in \ydac$ provided with a grading $M = \oplus_{j\in \Z} M^j$
such that
\begin{align*}
\Ac^i\cdot M^j &\subseteq M^{i+j}, & \delta(M^j) &\subseteq \sum_{h + k = j} \Ac^h \ot M^k, &i, j &\in \Z.
\end{align*}
The category $\grt\ydac$ of graded Yetter-Drinfeld modules over $\Ac$, with maps preserving all structures, is a
braided tensor category: the unit object $\ku$ has degree 0;
if $M, N \in\grt\ydac$, then $M \ot N \in\grt\ydac$ with the grading 
$(M \ot N)^j =\sum_{h + k = j} M^h \ot N^k$, and the braiding $c_{M,N}$ is homogeneous.

For instance, if $\Ac = \oplus_{i\in \Z} \Ac^i$ is finite dimensional, then $\Ac^* = \oplus_{i\in \Z} (\Ac^*)^i$ is also
a graded Hopf algebra, where $(\Ac^*)^i = (\Ac^{-i})^*$,  up to natural identifications. Then $D(\Ac)$ is also
a graded Hopf algebra, and $\grt\ydac$ is equivalent, as braided tensor category, to that of graded $D(\Ac)$-modules.

\begin{obs}\label{rem:coradical-graded}
	Let $\Ac = \oplus_{i\in \N_0} \Ac^i$ be a graded Hopf algebra. Then the coradical $\Ac_0$ coincides with the coradical of $\Ac^0$. In particular, if $\Ac^0$ is pointed, then so is $\Ac$ and  $G(\Ac) = G(\Ac^0)$.
\end{obs}

\pf The family $\Fc_n \Ac = \oplus_{i\in \I_{0,n}} \Ac^i$
is a coalgebra filtration  of $\Ac$. 
By \cite[Lemma 5.3.4]{Mo-libro}, $\Ac_0 \subseteq \Fc_0 \Ac =\Ac^0$, so $\Ac_0 \subseteq (\Ac^0)_0 \subseteq \Ac_0$.  
\epf

\medbreak
If $M \in\grt\ydac$, then $T(M)$ is a Hopf algebra in $\grt\ydac$, and so is $\toba(M)$ (because the quantum symmetrizer is homogeneous). 
We consider the $\Z^2$-gradings on these algebras given by
\begin{align*}
\deg T^n(M)^j &= \deg \toba^n(M)^j = (j, n),&  j&\in \Z, n \in \N_0.
\end{align*}	
Then $T(M)$ and $\toba(M)$ are $\Z^2$-graded algebras and coalgebras, but beware they are not 
\emph{$\Z^2$-graded algebras in} $\ydac$ as  $T^n(M)^j$ is not a Yetter-Drinfeld submodule. 
However we have:

\begin{lema}\label{lema:double-grading}
	The bosonizations $T(M) \# \Ac$ and $\toba(M) \# \Ac$ are $\Z^2$-graded Hopf algebras, 
	with grading given by
	\begin{align*}
	\deg \toba^n(M)^j \# \Ac^i &= (i+j, n),&  i,j&\in \Z, n \in \N_0.
	\end{align*}
\end{lema}

\pf This is straightforward:
\begin{align*}
&(\toba^n(M)^j \# \Ac^i)(\toba^m(M)^r \# \Ac^s)  \subseteq \sum_{h + k = i} \toba^n(M)^j \Ac^h \cdot \toba^m(M)^r \# \Ac^k \Ac^s \\ 
& \subseteq \sum_{h + k = i}  \toba^n(M)^j \toba^m(M)^{r+h} \# \Ac^{k+s} 
\subseteq \sum_{h + k = i}  \toba^{n+m}(M)^{j+r+h} \# \Ac^{k+s};
\\
&\Delta(\toba^n(M)^j\# \Ac^i) \subseteq \sum_{\substack{ p + q = n \\ h + k = i \\ r+s+t = j}}
\toba^p(M)^r  \#  \Ac^s   \Ac^h \ot \toba^q(M)^t \#    \Ac^k.
\end{align*}
\epf

\subsubsection{Filtered}\label{subsubsec:filtered-yd}
As above, let $\Ac = \oplus_{i\in \Z} \Ac^i$ be a graded Hopf algebra.
A \emph{filtered} Yetter-Drinfeld module over $\Ac$ is $M \in \ydac$ provided with an ascending filtration $(M\_{j})_{j\in \Z}$
such that $M = \bigcup_{j\in \Z} M\_{j}$ and
\begin{align*}
\Ac^i\cdot M\_{j} &\subseteq M\_{i+j}, & \delta(M\_{j}) &\subseteq \sum_{h + k = j} \Ac^h \ot M\_{k}, &i, j &\in \Z.
\end{align*}
The category $\filt\ydac$ of filtered Yetter-Drinfeld modules over $\Ac$, with maps preserving all structures, is a
braided tensor category: if $M, N \in\filt\ydac$, then $M \ot N \in\filt\ydac$ with filtration 
\begin{align}\label{eq:filtration-tensor-prod}
(M \ot N)\_{j} &=\sum_{h + k = j} M\_h \ot N\_k, 
\end{align}
so that the braiding $c_{M,N}$ preserves the filtration.

\begin{obs}\label{rem:filtered-graded}
	\begin{enumerate}[leftmargin=*,label=\rm{(\alph*)}]
		\item\label{item:graded-filtered} There is a forgetful functor $\grt\ydac\to \filt\ydac$:
		\begin{align*}
M&=\oplus_{k\in \Z} M^k \in \grt\ydac \rightsquigarrow &
&(M\_{j})_{j\in \Z}, & M\_{j} &=\oplus_{k \le j} M^k, &  j &\in\Z.
		\end{align*}
		\item\label{item:filtered-projection} If $\pi:M\to N$, is a projection in $\ydac$ and $M = \bigcup_{j\in \Z} M\_{j}$ is filtered, then $N$ is also filtered: $N = \bigcup_{j\in \Z} N\_{j}$, $N\_{j}:=\pi(M\_{j})$.
		\item\label{item:filtered-quotient-Verma} As a consequence, every quotient of a Verma module $M(\lambda)$ as in Proposition \ref{prop:L(lambda)} is a filtered Yetter-Drinfeld module.
	\end{enumerate}
\end{obs}

\medbreak
If $M \in\filt\ydac$, then $T^n(M)=M^{\otimes n} \in \filt\ydac$ by \eqref{eq:filtration-tensor-prod}, and also $\toba^n(M) \in \filt \ydac$ by Remark \ref{rem:filtered-graded} \ref{item:filtered-projection}. 
By \eqref{eq:filtration-tensor-prod} again,
\begin{align*}
T^m(M)\_i T^n(M)\_{j} \subseteq T^{m+n}(M)\_{i+j},
\end{align*}
and since the braiding preserves the filtration:
\begin{align*}
\Delta(T^m(M)\_i) \subseteq \sum_{\substack{n+p=m \\ j+k=i}} T^n(M)\_{j} \ot T^p(M)\_k.
\end{align*}
Thus the corresponding properties hold in $\toba(V)$ by Remark \ref{rem:filtered-graded} \ref{item:filtered-projection}:
\begin{align}\label{eq:B(M)-prod-coprod}
\begin{aligned}
\toba^m(M)\_i \toba^n(M)\_{j} & \subseteq \toba^{m+n}(M)\_{i+j},
\\
\Delta(\toba^m(M)\_i) & \subseteq \sum_{\substack{n+p=m \\ j+k=i}} \toba^n(M)\_{j} \ot \toba^p(M)\_k.
\end{aligned}
\end{align}

\begin{lemma}\label{lemma:bosonization-filtered}
	The Hopf algebra $\toba(M) \# \Ac$ is filtered, 
	with filtration
	\begin{align}\label{eq:bosonization-filtered}
	(\toba(M) \# \Ac)\_n &= \sum_{i+j+k\le n} \toba^i(M)\_{j} \# \Ac^k,&  n \in \Z.
	\end{align}
\end{lemma}

\pf By \eqref{eq:B(M)-prod-coprod} the filtration \eqref{eq:bosonization-filtered} is of algebras and coalgebras:
\begin{align*}
&(\toba^n(M)\_{j} \# \Ac^i)(\toba^m(M)\_r \# \Ac^s)  
\subseteq \sum_{h + k = i} \toba^n(M)\_{j} \Ac^h \cdot \toba^m(M)\_r \# \Ac^k \Ac^s 
\\ 
& \subseteq \sum_{h + k = i}  \toba^n(M)\_{j} \toba^m(M)\_{r+h} \# \Ac^{k+s} 
\subseteq \sum_{h + k = i}  \toba^{n+m}(M)\_{j+r+h} \# \Ac^{k+s};
\\
&\Delta(\toba^n(M)\_{j}\# \Ac^i) \subseteq \sum_{\substack{ p + q = n \\ h + k = i \\ r+s+t = j}}
\toba^p(M)\_r  \#  \Ac^s   \Ac^h \ot \toba^q(M)\_t \#    \Ac^k.
\end{align*}
Finally the antipode preserves the filtration since
\begin{align*}
\Ss (\toba^n(M)\_{j}\# \Ac^i) & \subseteq \sum_{h + k = j} (1\# \Ss_h(\Ac^h \Ac^i)) (\Ss_{\toba(M)}(\toba^n(M)\_k) \# 1)
\\
& \subseteq \sum_{h + k = j} (1\# \Ac^{h+i}) (\toba^n(M)\_k \# 1)
\subseteq (\toba(M)\# \Ac)\_{n+j+i}.
\end{align*}
That is,  \eqref{eq:bosonization-filtered} is a Hopf algebra filtration.
\epf

\begin{obs}\label{rem:bosonization-prenichols-filtered}
The statement holds also for pre-Nichols algebras of $M$, i.e. graded connected Hopf algebras in $\ydav$
generated by the degree one component which is isomorphic to $M$.
\end{obs}

\subsection{Pointed Hopf algebras}\label{subsec:hopf-subalg}
We state two facts needed later in the proof of the main result.

Let $B$ be a finite-dimensional pointed Hopf algebra with  $G = G(B)$ abelian. 
The graded Hopf algebra $\gr B = \oplus_{n\ge 0} \gr^n B$ associated to the coradical filtration of $B$ 
splits as $\gr B \simeq R \# \ku G$; $R = \oplus_{n\ge 0} R^n$ is a graded Hopf algebra in $\ydG$ and 
$\bV = R^1$ is the infinitesimal braiding of $B$. 
Now the projection $B_1 \to B_1/B_0$ admits a section $\mathfrak s: B_1/B_0 \to B_1$ of $\ku G$-bimodules.
In fact
\begin{align}\label{eq:first-term-splitting}
 B_1 = \ku G \oplus \mathfrak s(\bV) \ku G.
\end{align}

Let $\theta :=\dim \bV$, $\I = \I_{\theta}$. 
We fix a basis $(x_i)_{i\in \I}$ of $\bV$ such that 
the $G$-action and the $G$-coaction on $x_i$ are 
given by $\chi_i \in \widehat{G}$ and $g_i \in G$, respectively, for all $i\in \I$. 
This choice of the basis induces $\Z^\I$-gradings on $T(\bV)$ and $\toba(\bV)$. 
If $r\in T(\bV)$ is $\Z^\I$-homogeneous  of degree $(a_i)_{i\in\I}$, then we set $\chi_r=\prod_{i \in \I} \chi_i^{a_i}$.
Quite a bit is known about the
structure of  $B$:

\begin{itemize}[leftmargin=*]\renewcommand{\labelitemi}{$\circ$}
\item The list of all $\bV$ with $\dim \toba(\bV)< \infty$ is known by \cite{H-classif RS}.

\item By \cite{A-presentation} $R \simeq \toba(\bV)$
as mentioned in the Introduction. 

\item In \cite{A-jems,A-presentation} there was described a minimal set of homogeneous relations $\strat\subset T(\bV)$ defining $\toba(\bV)$. 
	
\item 	By \cite[Theorem 1.1]{AG} (as a culmination of a series of papers with other authors), there exists a family $\bsl=(\lambda_r)_{r\in\strat}$  in $\Bbbk^{\strat}$ with the restriction
$\lambda_r=0$ if $\chi_r\neq \epsilon$, 
such that $K \simeq \uf(\bsl)$; here $\uf(\bsl)$ is presented as the quotient $T(\bV) \# \ku G$ by relations 
\begin{align*}
r & - p_r(\bsl), & r\in \strat
\end{align*}
where $p_r(\bsl) \in \oplus_{k < \deg r} T^k(\bV) \# \ku G$ are
defined recursively.
\end{itemize}

We shall use a refinement of these facts to prove the next result; we were unable to find a simpler proof.

\begin{lema}\label{lema:dim-subalg} 
	Let $B$ be a finite-dimensional pointed Hopf algebra with  $G = G(B)$ abelian; keep the notation above.
Let $\bV' \leq \bV$ in $\ydG$ and let $B'$ be the subalgebra of $B$ generated by 
$\mathfrak s(\bV')$ and $G$. Then $B'$ is a lifting of $\toba(\bV')\# \Bbbk G$, particularly
\begin{align}\label{eq:dim-subalg}
\dim B' = \dim \toba(\bV') \vert G \vert.
\end{align}
\end{lema}

\pf Evidently $B'$ is a Hopf subalgebra of $B$, in particular it is pointed. 
Since $\bV' \hookrightarrow$ the infinitesimal braiding of $B'$, we have 
\begin{align}\label{eq:dim-subalg1}
\dim B' \geq \dim \toba(\bV') \vert G \vert.
\end{align}
Let $t:=\dim \bV'\leq  \theta$. By  \cite[Lemma 3.7]{AG-survey}, there exist
\begin{itemize}[leftmargin=*]\renewcommand{\labelitemi}{$\circ$}
	\item a basis $(x_i)_{i\in \I}$ of $\bV$ as described above, such that $(x_i)_{i\in \I_t}$ is a basis of $\bV'$;
	
	\item a minimal set of homogeneous (with respect to the basis above) relations 
	$\strat\subset T(\bV)$ defining $\toba(\bV)$, such that $\strat' =\strat \cap T(V')$, (i.e. the relations with only letters $(x_i)_{i\in\I_{t}}$) is a minimal set of relations defining $\toba(V')$;
	
	\item and such that if $\bsl=(\lambda_r)_{r\in\strat}$ is a family in $\Bbbk^{\strat}$ as above
	such that $B \simeq \uf(\bsl)$, then $\bsl':=(\lambda_r)_{r\in\strat'}$ satisfies 
\begin{align*}
p_r(\bsl) & = p_r(\bsl'), & r\in \strat'
\end{align*}
and $\uf(\bsl') = T(\bV') \# \ku G / \langle r - p_r(\bsl'): r\in \strat'\rangle$ is a a lifting 
of $\toba(\bV')\# \Bbbk G$. 
\end{itemize}
Since clearly there is a surjective map
$\uf(\bsl')  \to B'$, we  have
\begin{align}\label{eq:dim-subalg2}
\dim B' \leq \dim \uf(\bsl') = \dim \toba(\bV') \vert G \vert
\end{align}
that together with \eqref{eq:dim-subalg1} gives \eqref{eq:dim-subalg}.
\epf

\begin{lema}\label{lema:primitivos-infbraiding}
	Let $T = \oplus_{n\ge 0} T^n$ be a finite-dimensional graded connected Hopf algebra in $\ydG$
	and $B = T \# \ku G$, so that $B$ is pointed. 
	Then there is a monomorphism $\Pc(T) \hookrightarrow \bV$, the infinitesimal braiding of $B$.
\end{lema}

\pf
Clearly there is a monomorphism $\Pc(T) \hookrightarrow B_1$; then apply \eqref{eq:first-term-splitting}.
\epf

\subsection{Semisimplicity}\label{subsec:pf-thm.Z-semisimple}

Again we assume that $H$ and $V$ satisfy \eqref{eq:hypothesis-1} and \eqref{eq:hypothesis-2}.
We consider the following question:

\begin{question}\label{question:fd-semisimple}
	Let $Z \in \ydhh$. Does the finiteness of 
	$\dim \toba(Z)$  imply that $Z$ is semisimple?
\end{question} 

\begin{obs}\label{rem:simple-quotient} If $Z$ provides a negative answer to Question \ref{question:fd-semisimple}, then we may assume that $Z/S$ is simple, where $S$ is the socle of $Z$. 
\end{obs}

Indeed, let $L$ be a simple subobject of $Z/S$ and consider $\widetilde{Z}=p^{-1}(L)\subset Z$, where $p:Z\to Z/S$ is the projection. Then the socle of $\widetilde{Z}$ is again $S$, $\toba(\widetilde{Z}) \subset \toba(Z)$ is finite-dimensional and $\widetilde{Z}/S\simeq L$, so we replace $Z$ by $\widetilde{Z}$.

\medbreak
We answer affirmatively Question \ref{question:fd-semisimple} when $H$ is the group algebra 
of a finite abelian group.

\begin{theorem}\label{thm:Z-semisimple}
Let $H = \ku G$, where $G$ is a finite abelian group, and let $V$ satisfy \eqref{eq:hypothesis-2}.  
Set $\Ac(V) = \toba(V) \# \ku G$.
Let $Z \in \ydav$. If $\dim \toba (Z) < \infty$, then $Z$ is semisimple.
\end{theorem}

We start by a general Lemma.

\begin{lemma}\label{lemma:primitivos-bosonization}
Let $H$ be an arbitrary  Hopf algebra, $R$  a Hopf algebra in $\ydhh$ and $\Hc = R\# H$. Then
$\Pc_{g,1}(\Hc)\subset \Pc(R)\#1 + 1\# \Pc_{g,1}(H)$ for all $g\in G$.
\end{lemma}
\pf
We fix a basis $(x_i)_{i \in I}$ of $R$ such that for a distinguished $i_0 \in I$, $x_{i_0} = 1$
and $\epsilon(x_i) = 0$ when $i \neq i_0$.
Let $w=\sum_{i\in I} x_i\# k_i \in \Pc_{g,1}(\Hc)$, where $k_i\in H$, $i\in I$.  
We apply $\id\#\id\ot\epsilon\#\id$ to the equality 
\begin{align*}\Delta(w) &= \sum_{i\in I} (x_i)^{(1)} \# (x_i)^{(2)}_{(-1)} (k_i)_{(1)} \otimes
(x_i)^{(2)}_{(0)}\# (k_i)_{(2)}
= w\otimes 1\#1 + 1\# g \otimes w\end{align*} and get
\begin{align*}
\sum_{i\in I} x_i\#  (k_i)_{(1)} \otimes  (k_i)_{(2)}
= 1 \# k_{i_0} \otimes 1 + \sum_{i_0 \neq i\in I} x_i\# k_i \otimes 1
+ 1\# g \otimes  k_{i_0},
\end{align*}
hence 
\begin{align*}
\Delta (k_i) &= k_i \otimes 1, \ i\in I, \ i\neq i_0; &
\Delta (k_{i_0}) &= k_{i_0} \otimes 1+  g \otimes  k_{i_0}.
\end{align*}
That is, $k_i \in \ku$ for $i\neq i_0$; $k_{i_0} \in \Pc_{g,1}(H) \subset \Pc_{g,1}(\Hc)$
and $w = z \# 1 + 1\# k_{i_0}$, where $z = \sum_{i_0 \neq i\in I} x_i  k_i$.
Then $z \#1 \in\Pc_{g,1}(\Hc)$. 
We apply $\id\#\epsilon\ot\id\#\id$ to the equality 
\begin{align*}\Delta(z \#1) &=  (z)^{(1)} \# (z)^{(2)}_{(-1)}  \otimes
(z)^{(2)}_{(0)}\# 1 = z \#1 \otimes 1\#1 + 1\# g \otimes z\#1 \end{align*} and get
\begin{align*}\Delta(z) &=  (z)^{(1)}  \otimes
(z)^{(2)} = z \otimes 1 + 1 \otimes z. \end{align*}
Thus $z \in \Pc(R)$.
\epf

\medbreak
Let $H$ and $V$ be as in \eqref{eq:hypothesis-1} and \eqref{eq:hypothesis-2}. Then 
$\Ac(V) = \oplus_{i \ge 0} \Ac^i(V)$ is a graded Hopf algebra where $\Ac^i(V) = \toba^i(V) \# H$.
Here is the key step.

\begin{lema}\label{lema:braided-bos-prenichols-findim}
	Assume that $H = \ku G$, where $G$ is a finite abelian group.
	Let $N \subseteq M_+(\lambda)$ be  a submodule of $M(\lambda)$ 
	and let $M = M(\lambda) / N \in {}^{\Ac(V)}_{\Ac(V)}\mathcal{YD}$. 
	If $\dim \toba (M) < \infty$, then $ M = L(\lambda)$.
\end{lema}

Thus the braided bosonization of a Nichols algebra does not need to be a Nichols algebra.
Compare with \cite[Thm. 4.3.1]{U}.

\pf
Let $J = \toba(M)\# \Ac(V)$. 
The filtration of  $M(\lambda)$ arising from its grading induces a filtration of $M$, and thus a Hopf algebra filtration on $J$. By \eqref{eq:bosonization-filtered}, the first two terms of this filtration are
\begin{align}\label{eq:coradical-J}
J\_0 &= \Bbbk \# \Bbbk \# \ku G \simeq \ku G, &
J\_1 &= \Bbbk \# \Bbbk \# \ku G \oplus \ku \# V \# \ku G \oplus  \lambda \# \ku \# \ku G. 
\end{align}
Hence  $J_0=\Bbbk G$ by \cite[Lemma 5.3.4]{Mo-libro} and $J$ is pointed.  
Let $\bV \in\ydG$ be the infinitesimal braiding of $J$; by \eqref{eq:coradical-J}
there is a monomorphism  
\begin{align*}
\bV' := V\oplus \lambda \hookrightarrow \bV.
\end{align*}
Let $J'$ be the subalgebra of $J$ generated by $\bV'$ and $G$, that is by
$J_1$. We claim that $J' = J$. Indeed 
$J$ is generated by $ M, V, G$ but $M = \ad \toba(V) \cdot \lambda$, 
implying the claim. By Lemma \ref{lema:dim-subalg} and the last claim, 
\begin{align*}
\dim \toba(\bV') \vert G \vert = \dim J'= \dim J = \dim \toba(M) \dim \toba(V) |G|.
\end{align*}
By Proposition \ref{prop:B(L(lambda))}, $\dim \toba(\bV') =\dim \toba(V)\dim \toba(L(\lambda))$, 
implying that  $\dim \toba(L(\lambda)) = \dim \toba(M)$. As $M$ projects onto $L(\lambda)$, we conclude that $M=L(\lambda)$.
\epf

We are ready to complete the proof.

\subsubsection*{Proof of Theorem \ref{thm:Z-semisimple}.}
We fix some notation:

\medbreak
\begin{itemize}[leftmargin=*]\renewcommand{\labelitemi}{$\circ$}
\item $S$ is the socle of $Z$, $S=L(\lambda_1) \oplus \dots \oplus L(\lambda_t)$ for some $\lambda_i \in \Irr \ydG$.
	
\medbreak
\item $\Hc:= \toba(Z) \# \Ac(V)$, $\Hc':=\toba(S) \# \Ac(V)$. We may consider $\Hc'$ as a Hopf subalgebra of $\Hc$, since $\toba(S)$ is a Hopf subalgebra of $\toba(Z)$ in $\ydav$.
\end{itemize}
By Remark \ref{rem:coradical-graded}, $\Hc$ is pointed and $G(\Hc) = G$.
Thus $\Hc$ is generated by $\Hc_0=\ku G$ and the skew-primitive elements. 
Now we want to describe $\Pc_{g,1}(\Hc)$ for each $g\in G$: By Lemma \ref{lemma:primitivos-bosonization}, with $R = \toba (Z)$,
\begin{align*}
\Pc_{g,1}(\Hc)&\subset Z\#1 + 1\# \Pc_{g,1}(H),& g&\in G.
\end{align*}

\medskip
We next determine the elements of $Z$ that are skew-primitive in $\Hc$. Set
\begin{align*}
Z' &=\sum_{g\in G} Z'_g, & \text{where }Z'_g&:=\{z\in Z \, | \, z\#1 \in \Pc_{g,1}(\Hc) \}, \, g\in G.
\end{align*}
We check that
\begin{align*}
h\cdot z & \in Z'_{hgh^{-1}} = Z'_g, & \delta_h\cdot z&= \delta_{g,h} z, & \text{for all }g,h\in G, & \, z\in Z'_g.
\end{align*}
Hence $Z'$ is a $D(\Bbbk G)$-submodule of $Z$. Moreover, $Z' =\oplus_{g\in G} Z'_g$.

\begin{step}
	$F\cdot z=0$ for all $F\in \overline{V}$ and all $z\in Z'$.
\end{step}

Let $(E_k)_{k\in \I_M}$ be a basis of $\toba(V)$, where $E_1=1$, each $E_k$ is homogeneous with respect to the $\N_0$-grading and the $\ku G$-coaction is given by $g_k\in G$. Let $(F_k)_{k\in \I_M}$ be the dual basis of $\toba(\overline{V})$.
Let $z\in Z'_g$, $g\in G$. Hence,
\begin{align*}
\Delta(z\# 1) & \overset{\star}{=} z\# 1 \otimes 1\# 1 +  1\# g \otimes z\# 1
\\
& \overset{\ast}{=} z\# 1 \otimes 1\# 1 + \sum_{k\in \I_M} 1\# E_k gg_k^{-1} \otimes F_k\cdot z\# 1;
\end{align*}
here $\star$ is merely by assumption, while $\ast$ comes from the formulas for the bosonization.
Thus the Step follows.

\begin{step}
	$Z' \subset S$.
\end{step}

We decompose $Z'=\oplus_{i\in \Ib_{t}} \mu_i$, each $\mu_i$ a simple $D(\Bbbk G)$-module.
Let $z\in \mu_i$, $z\neq 0$, $i\in \Ib_{t}$. By the previous Step, $F\cdot z=0$ for all $F\in \overline{V}$. Thus, there exists a surjective map $M(\mu_i)\to \doble \cdot z=\doble^{\ge 0} \cdot z$. Also, $\toba(\doble \cdot z) \subseteq \toba(Z)$, so $\dim \toba(\doble \cdot z)<\infty$. 
By Lemma \ref{lema:braided-bos-prenichols-findim}, $\doble \cdot z \simeq L(\mu_i)$, so $z\in S$.

\begin{step}
	The subalgebra generated by $\Hc_0 \cup (\cup_{g\in G} \Pc_{g,1}(\Hc))$ is $\Hc'$.
\end{step}

From the previous steps, the skew-primitive elements of $\Hc$ belong either to $S\# 1$ or else to $1\# \Ac(V)$. Notice that
\begin{align*}
\Hc'&=\toba(S) \# \Ac(V) = \toba(S) \# (\toba(V) \# \ku G) \simeq \toba(V \oplus \lambda_1 \oplus \dots \oplus\lambda_t) \# \ku G.
\end{align*}
Thus the subalgebra generated by $\ku G$ and the skew-primitive elements is $\Hc'$. 

\medbreak

By the last Step, $\Hc'= \Hc$. As $\Hc'\cap Z\#1=S\#1$, we have that $S=Z$.
\qed

\subsection{The diagram}\label{subsec:diagram}

We consider the following setting:

\begin{itemize}[leftmargin=*]\renewcommand{\labelitemi}{$\circ$}
	
	\medbreak
	\item $G$ is a finite abelian group and $V \in \ydG$ satisfies  
 \eqref{eq:hypothesis-2}; $\Ac(V) = \toba(V) \# \ku G = \oplus_{i\in \N_0} \toba^i(V) \# \ku G$.
		
	\medbreak
	\item $\Rt = \oplus_{i\in \N_0} \Rt^i$ is a finite-dimensional connected graded Hopf algebra in $\ydav$, generated by 
$Z := \Rt^1$. Thus we have an epimorphism $\Rt\to \toba(Z)$ of Hopf algebras in $\ydav$. In other words, $\Rt$ is a pre-Nichols algebra of $Z$.
\end{itemize}

\begin{prop}\label{prop:R-nichols}
$\Rt\simeq \toba(Z)$.
\end{prop}

\pf By Theorem \ref{thm:Z-semisimple}, $Z$ is semisimple, say $Z = L(\lambda_1) \oplus \dots \oplus L(\lambda_t)$, 
where $\lambda_1, \dots, \lambda_t \in \Irr \ydG$. Particularly, $Z$ is a graded object of $\ydav$, hence so is $T(Z)$,
and consequently $\Rt$ is a filtered Hopf algebra in $\ydk$, and \emph{a fortiori} $J := \Rt \# \Ac(V)$ is a filtered Hopf algebra, with filtration given by
\begin{align}\label{eq:bosonization-filtered-Rt}
J_{(n)} &= \sum_{i+j+k\le n} \Rt^i\_{j} \# \toba^k(V) \# \ku G,&  n \in \N_0,
\end{align}
by (the same proof as) Lemma \ref{lemma:bosonization-filtered}, cf. Remark \ref{rem:bosonization-prenichols-filtered}. 
Now the proof follows as the one of Lemma \ref{lema:braided-bos-prenichols-findim}. 
Indeed $\Rt \# \Ac(V)$ is pointed by \cite[Lemma 5.3.4]{Mo-libro}, with coradical
$\Bbbk G$.
Let $\bV$ be the infinitesimal braiding of $J$. Then
$$\bV' := V\oplus \lambda_1 \oplus \dots \oplus \lambda_t \hookrightarrow \bV.$$ 

Let $J'$ be the subalgebra of $J$ generated by $\bV'$ and $G$. We claim that $J' = J$. Indeed 
$J$ is generated by $ Z, V, G$ but $Z = \ad \toba(V) \cdot (\lambda_1 \oplus \dots \oplus \lambda_t)$, 
implying the claim. By Lemma \ref{lema:dim-subalg} and the last claim, 
\begin{align*}
\dim \toba(\bV') \vert G \vert = \dim J'= \dim J = \dim \Rt \dim \toba(V) |G|.
\end{align*}
By Proposition \ref{prop:B(L(lambda))}, $\dim \toba(\bV')=\dim \toba(V)\dim \toba(Z)$, so $\Rt \simeq \toba(Z)$.
\epf

We need one more general result before proving the main Theorem.

\begin{lemma}\label{lemma:prim-standard-filtration}
Let $\hopf$ be a Hopf algebra with Hopf coradical $\hopf_{[0]} \simeq H$, so that $\gr \hopf \simeq R \# H$ for some connected graded Hopf algebra $R = \oplus_{i\in \N_0} R^i \in \ydhh$. Then $\Pc(R)=R^1$.
\end{lemma}
\pf
As $R$ is connected and graded, $R^1\subseteq \Pc(R)$. On the other hand, let $x\in\Pc(R)$. By \cite[Proposition 1.2]{AC},
\begin{align*}
(\gr \hopf)_{[0]} &= R^0 \# H \simeq H, &
(\gr \hopf)_{[1]} &= R^0 \# H \oplus R^1 \# H.
\end{align*}
By direct computation,
\begin{align*}
\Delta(x\#1) &= x\#1  \otimes 1\#1 + 1\# x_{(-1)}  \otimes x_{(0)}\#1
\\
& \in \gr \hopf \otimes (\gr \hopf)_{[0]} 
+ (\gr \hopf)_{[0]} \otimes \gr \hopf.
\end{align*}
Hence $x\# 1 \in (\gr \hopf)_{[0]} \wedge (\gr \hopf)_{[0]} 
= (\gr \hopf)_{[1]}$, so $x\in R^1$.
\epf

If $H$ is cosemisimple, then it follows from the previous result that $R$ is coradically graded \cite{AS-98}.
But we do not know whether this is true in general.

\subsubsection*{Proof of Theorem \ref{thm:main}}

Recall that $L$ is a basic finite-dimensional Hopf algebra such that $G = \Hom_{\text{alg}} (L, \ku)$ is abelian
and  $\hopf$ is a Hopf algebra with $\hopf_{[0]} \simeq L$, so that $\gr \hopf \simeq R \# L$.
Also $B = L^*$ and $K = \gr B = \toba(V) \# \ku G$. Since \eqref{item:intro-basic-hyp2} holds, 
we have inverse equivalences of braided tensor categories $\Fc: \ydl \to \ydk$ and
$\Gc: \ydk \to \ydl$ that preserve dimensions, gradings, Hopf algebras, spaces of primitive elements
and Nichols algebras. 

\medbreak
The implication \ref{item:intro-main-2} $\implies$ \ref{item:intro-main-1} is a consequence of Proposition \ref{prop:B(L(lambda))}.

\medbreak
\ref{item:intro-main-1} $\implies$ \ref{item:intro-main-2}: 
we consider the graded Hopf algebra $\Fc(R) =\oplus_{i\in \N_0} \Fc(R)^i \in \ydk$, where $\Fc(R)^i=\Fc(R^i)$. 
Since $\Pc(\Fc(R)) = \Fc(\Pc(R))$,  Lemma \ref{lemma:prim-standard-filtration} says that $\Pc(\Fc(R)) =\Fc(R)^1 =: Z$. 
Let $\Rt$ be the subalgebra of $\Fc(R)$ generated by $Z$. 
By Theorem \ref{thm:Z-semisimple}, $Z$ is semisimple. Hence
$Z = L(\lambda_1) \oplus \dots \oplus L(\lambda_t)$ for some $\lambda_1, \dots, \lambda_t \in \Irr \ydG$ and \eqref{eq:intro-main2} holds by Theorem \ref{theorem:intro-basic}.

By Proposition \ref{prop:R-nichols}, $\Rt\simeq \toba(Z)$. 
By Remark \ref{rem:coradical-graded}, $\Fc(R)\# K$ is pointed; being finite-dimensional, 
it is generated by $G$ and the subspaces $\Pc_{g,1}$, $g\in G$. By Lemma \ref{lemma:primitivos-bosonization},
$\Pc_{g,1}\subset \Pc(\Fc(R))\# 1 + 1\# V$ for all $g\in G$. 
As $\Pc(\Fc(R))=Z$, $\Fc(R)\# K$ is generated by $Z$, $V$ and $G$; as
$\ku\langle Z, V, G \rangle = \Rt\# K$, we have that $\Fc(R)\# K = \Rt\# K$.
That is, $\Fc(R) \simeq \toba(Z)$. Therefore
\begin{align*}
R & \simeq \Gc(\Fc(R)) \simeq \Gc(\toba(Z)) \simeq \toba (\Zc),
& \Zc &= \Gc(L(\lambda_1)) \oplus \dots \oplus \Gc(L(\lambda_t)),
\end{align*}
which completes the proof.
\qed

\section{Decompositions from diagonal type}\label{sec:nichols-diag}
In this Section we compute some of the Nichols algebras $\toba(Z_U)$ assuming that $W = V \oplus U$ 
is  of diagonal  type.

\subsection{Generalities}\label{subsec:generalities}
\subsubsection{The setting}\label{subsubsec:generalities-setting}
Let $\theta \in \N$ and $\I = \I_{\theta}$.
Let $W$ be a braided vector space of diagonal type, of dimension $\theta$ with braiding matrix $(q_{ij})_{i, j \in \I}$
with respect to a basis $(x_{i})_{i \in \I}$. This matrix 
defines a $\Z$-bilinear form $\bq:\Z^{\I}\times\Z^{\I}\to\ku^\times$ by $\bq(\alpha_j,\alpha_k)=q_{jk}$ for all $j,k\in\I$.
We set 
\begin{align*}
\bq_{\alpha \beta}  &= \bq(\alpha,\beta),& N_{\beta} &= \ord \bq_{\beta \beta} \in \N \cup \{\infty\},&
\alpha,\beta  &\in \Z^{\I},\ i \in \I.
\end{align*}

We assume that $q_{ii} \neq 1$, $i\in \I$ and set $\widetilde{q}_{ij} = q_{ij}q_{ji}$, $i< j$.
We keep the notations from \cite{AA-diag-survey}.
We describe the braiding of $W$ by its Dynkin diagram \cite{H-classif RS}, see e.g. \cite[\S 2.1]{AA-diag-survey}; for instance, if $\theta = 2$, 
$\Dynkintwoxy{q_{11}}{1}{\widetilde{q}_{12}}{q_{22}}{2}$. 

Let $\Gamma$ be a free abelian group of rank $\theta$ with a fixed basis $(g_i)_{i\in \I_\theta}$.
We define $\chi_i: \Gamma \to \ku^{\times}$ by 
$\chi_i(g_j) = q_{ji}$. We realize $W \in \ydg$ by declaring that 
\begin{align*}
x_i &\in W_{g_i}^{\chi_i},& i &\in \I_\theta.
\end{align*}
As in \cite{AA-diag-survey}, we set 
\begin{align}\label{eq:xij}
x_{ij} &= \ad_c x_i (x_j) = x_ix_j - q_{ij} x_jx_i,& i&\neq j \in \I,
\\ \label{eq:iterated}
x_{i_1i_2\cdots i_k}& =(\ad_c x_{i_1})\cdots(\ad_c x_{i_{k-1}})\, (x_{i_k}),& &i_1, i_2, \cdots, i_k\in\I,
\\ \label{eq:roots-Atheta}
x_{(k \, l)} &= x_{k\,(k+1)\, (k+2) \dots l},& &k < l \in \I.
\end{align}

We fix $t \in \I_{\theta - 1}$; let $V$, respectively $U$, be the subspace generated by $(x_{i})_{i \in \I_t}$,
respectively $(x_{i})_{i \in \I_{t+1, \theta}}$.
We keep the notation in \eqref{eq:braided-bosonization}.
Then 
\begin{align*}
\toba(W) \simeq \toba (Z_U) \# \toba (V)\simeq \toba (Z_V) \# \toba (U).
\end{align*}

\subsubsection{Assumptions}\label{subsubsec:generalities-assumptions}
We assume that $W$ is arithmetic, i.e. that the set $\varDelta_+$ of positive roots of $W$ is finite, see \cite[\S 2.4]{AA-diag-survey}.
Then the $n_{ij}$'s defined below are integers, i.e. the sets on the right are non-empty:
\begin{align*}
n_{ij}&:=  \min \left\{ n \in \N_0: (n+1)_{q_{ii}}
(1-q_{ii}^n q_{ij}q_{ji} )=0 \right\},& i&\neq j \in \I.
\end{align*}

We fix a convex order of $\varDelta_+$ \cite{A-jems}, see \cite[Theorem 2.36]{AA-diag-survey}. Correspondingly, 
we have the numeration   $(\beta_k)_{k \in \I_{\ell}}$ of  $\varDelta_+$.
For every $k \in \I_{\ell}$, let $x_{\beta_k}$ be the corresponding root vector; these are defined as in \cite[Remark 2.14]{AA-diag-survey}
when the convex order arises from Lyndon words, or via the Lusztig isomorphisms in general.
These generate a PBW-basis of $\toba(W)$, that is
\begin{align}\label{eq:pbw}
\left\{ x_{\beta_{\ell}}^{n_{\ell}} x_{\beta_{\ell - 1}}^{n_{\ell - 1}} \dots x_{\beta_2}^{n_{2}}  x_{\beta_1}^{n_{1}} \, | \,   
n_{k} \in \I_{0, N_{\beta_k} - 1}, k \in \I_{\ell} \right\}.
\end{align}

Let $i<j \in \I_{\ell}$, $p_{i+1}, \dots, p_{j-1} \in\N_0$. Since the order is convex, 
there exist $c_{p_{i+1}, \dots, p_{j-1}}^{(i,j)} \in \ku$, explicitly computable \cite[Lemma 4.5]{A-jems}, such that
\begin{equation}\label{eq:quantumSerregeneralizadas}
\left[ x_{\beta_i}, x_{\beta_j} \right]_c= 
\sum_{p_{i+1}, \dots, p_{j-1} \in\N_0} c_{p_{i+1}, \dots, p_{j-1}}^{(i,j)} \ x_{\beta_{j-1}}^{p_{j-1}} \dots x_{\beta_{i+1}}^{p_{i+1}}.
\end{equation}
Notice that if $\sum p_k\beta_k \neq \beta_i+\beta_j$, 
then $c_{p_{i+1}, \dots, p_{j-1}}^{(i,j)}  =0$, since $\toba_{\bq}$ is $\N_0^{\I}$-graded.

\begin{obs}\label{rem:consecutive-roots}
Let $j = i + 1 \in \I_{\ell}$. Then $\left[ x_{\beta_i}, x_{\beta_j} \right]_c= 0$ by \eqref{eq:quantumSerregeneralizadas}.
\end{obs}

\subsubsection{Properties of $\toba(Z_U)$}\label{subsubsec:generalities-properties}
We first describe a basis of $\toba(Z_U)$. Given $\alpha = \sum_{i\in \I} c_i \alpha_i\in \Z^{\I}$, 
let $\supp \alpha = \{i\in \I: c_i \neq 0\}$. 
Let $\Delta_+^V$ be the set of positive roots of $V$, $\lgot=|\Delta_+^V|$. Hence $\Delta_+^V= \{ \alpha\in\Delta_+ : \supp \alpha\subset \I_t \}$.

\begin{lema}\label{lema:convex-order}
There exists a convex order $<$ on $\Delta_+$ such that
$\Delta_+^V=\{ \beta_k\}_{k\in\I_{\lgot}}$.
\end{lema}
\pf Here we shall use the notation of \cite[\S 2.7]{AA-diag-survey}; in particular, the Weyl groupoid of $V$ is denoted $\mathcal W_V$.
Let $\omega_0^V$ be the element of maximal length of $\mathcal W_V$ that ends at $V$. 
We fix a reduced expression $\omega_0^V=\sigma_{i_1}^V\sigma_{i_2} \dots \sigma_{i_\lgot}$, $i_j\in\I_t$. 
By abuse of notation, we consider the element $\omega = \sigma_{i_1}^W \sigma_{i_2} \dots \sigma_{i_\lgot}$ of $\mathcal W_W$ that ends at $W$.
Clearly, $\omega$ has length $\lgot$,
that is, the expression $\sigma_{i_1}^W \sigma_{i_2} \dots \sigma_{i_\lgot}$ since $s_{i_1}^W \dots s_{i_{j-1}} (\alpha_{i_j}) = s_{i_1}^V \dots s_{i_{j-1}} (\alpha_{i_j})$ is a positive root for all $j\in\I_{\lgot}$, cf. \cite[Corollary 3]{HY-Weyl-gpd}. Using \cite[Corollary 4]{HY-Weyl-gpd} we may extend this expression to an expression of the element of maximal length $\omega_0^W=\sigma_{i_1}^W \dots \sigma_{i_\lgot} \sigma_{i_{\lgot+1}} \dots \sigma_{i_\ell}$; hence the first $\lgot$ roots of the associated convex order are the roots in $\Delta_+^V$.
\epf


By a combination of the main results of \cite{U1,U2}, the algebra $\toba(Z_U)$ has a PBW-basis. 
We next give an explicit choice of such a basis.

\begin{prop}
The elements $x_{\beta_k}$ with  $k \in \I_{\lgot+1,\ell}$ generate a PBW-basis of $\toba(Z_U)$. That is,
the following set is a basis of $\toba(Z_U)$:
\begin{align}\label{eq:pbw-Z_U}
\left\{ x_{\beta_{\ell}}^{n_{\ell}} x_{\beta_{\ell - 1}}^{n_{\ell - 1}} \dots x_{\beta_{\lgot+1}}^{n_{\lgot+1}} \, | \,   
n_{k} \in \I_{0, N_{\beta_k} - 1}, k \in \I_{\lgot+1,\ell} \right\}.
\end{align}
The order of $\I_{\lgot+1,\ell}$ defines a filtration on $\toba(Z_U)$ whose associated graded algebra is a (truncated) quantum polynomial ring.
\end{prop}

\pf 
Let $i\in\I_t$, $j \in \I_{\lgot+1,\ell}$. Suppose that $\partial_i(x_{\beta_j})\neq 0$: that is, some $x_{\beta_{\ell}}^{m_{\ell}} x_{\beta_{\ell - 1}}^{m_{\ell - 1}} \dots x_{\beta_{1}}^{m_{1}} \ot x_i$ appears with non-zero coefficient in the expression of $\Delta(x_{\beta_j})$.
The subspace spanned by $\left\{ x_{\beta_{j}}^{n_{j}} x_{\beta_{j- 1}}^{n_{j- 1}} \dots x_{\beta_{1}}^{n_{1}} \, | \,   n_{k} \in \I_{0, N_{\beta_k}}, k \in \I_{j} \right\}$ is a right coideal subalgebra \cite{HS-coideal subalg}, so $m_{k}=0$ for all $k\in\I_{j+1,\ell}$; as the coproduct preserves the $\Z^{\I}$-degree, we also have $m_j=0$, and 
\begin{align*}
\beta_j=\alpha_i+ \sum_{k\in\I_{j-1}} m_k\beta_k.
\end{align*}
Note that $\alpha_i=\beta_k$ for some $k\in\I_{j-1}$, see Lemma \ref{lema:convex-order}. As the order is strongly convex \cite[Theorem 3.11]{A-jems}, we have that $\beta_1< \beta_j<\beta_{j-1}$, a contradiction. Hence $\partial_i(x_{\beta_j})=0$ for all $i\in\I_t$, and by Proposition \ref{prop:B(ZU)-derivations}, $x_{\beta_j} \in \toba(Z_U)$.

As $\toba(Z_U)$ is a subalgebra of $\toba(W)$, every element in \eqref{eq:pbw-Z_U} belongs to $\toba(Z_U)$. 
The set \eqref{eq:pbw-Z_U} is linearly independent and spans a  subspace with the same Hilbert series as $\toba(Z_U)$, hence it is a basis of $\toba(Z_U)$.

The last statement follows e.g. as \cite[Proposition 16, Corollary 17]{A-pre-Nichols}.
\epf

A first immediate consequence is that $\toba(Z_U)$ is noetherian.
If $W$ is generic, hence of Cartan type \cite{AA,R}, then it satisfies further properties.

\begin{theorem}\label{thm:ZU} Assume that $W$ is generic.
The algebra $\toba(Z_U)$ is an iterated Ore extension; thus it is strongly noetherian, AS-regular and Cohen-Macaulay.
\end{theorem}

In many cases the Theorem was already known, e.g. \cite[(0.2), (8.5)]{ArS}, \cite[Remark 2.9]{NS}, \cite[Theorem 4.3]{E}, \cite[\S 7]{LW}.

\pf First, we check that $\toba(Z_U)$ is an iterated Ore extension. In fact, for each $j\in\{\lgot+1,\ell\}$ let $R_j$ be the subspace spanned by
\begin{align*}
\left\{ x_{\beta_{\ell}}^{n_{\ell}} x_{\beta_{\ell - 1}}^{n_{\ell - 1}} \dots x_{\beta_{j}}^{n_{j}} \, | \,   n_{k} \in \I_{0, N_{\beta_k}}, k \in \I_{\lgot+1,\ell} \right\}.
\end{align*}
Thus each $R_{j}$ is a free $R_{j+1}$-module with basis $(x_{\beta_{j}}^n)_{n\in\N_0}$.
Let $\sigma_j$ be the algebra automorphism given by the action of $g_j$, and $D_j:= [x_{\beta_{j}}, -]_c$. Then $D_j$ is a $(\id, \sigma_j)$-derivation, and by \eqref{eq:quantumSerregeneralizadas}, $D_j(x_{\beta_k})\in R_{j+1}$ for all $k\in\I_{j+1,\ell}$. Hence 
\begin{align*}
[x_{\beta_{j}}, x]_c &=x_{\beta_{j}} x- (g_j\cdot x) x_{\beta_{j}} \in R_{j+1} & \text{for all }& x\in R_{j+1},
\end{align*}
and $D_j$ restricts to a $(\id, \sigma_j)$-derivation of $R_{j+1}$. Therefore 
\begin{align*}
R_j & \simeq R_{j+1}[x_{\beta_j};\sigma_j, D_j], & 
\text{for all } & j\in \I_{\lgot+1,\ell-1}.
\end{align*}

Now $\toba(Z_U)$ is strongly noetherian
by \cite[Proposition 4.10]{ASZ};  AS-regular by \cite[Proposition 2]{AST} and Cohen-Macaulay by \cite[Lemma 5.3]{ZZ}.
\epf

Next we deal with the Laistrygonian Nichols algebras $\toba(\lstr( 1, \ghost))$ \cite{AAH}.

\begin{theorem}\label{thm:laistrygonian}
	The algebra $\toba(\lstr( 1, \ghost))$ is an iterated Ore extension; thus it is strongly noetherian, AS-regular and Cohen-Macaulay.
\end{theorem}
\pf 
We follow the notation of \cite[\S 4.3]{AAH}. The subalgebra of $\toba(\lstr( 1, \ghost))$ generated by $x_1$, $x_2$ is (isomorphic to) the Jordan plane, so it is an Ore extension of the polynomial algebra in one variable. 

For each $j\in\I_{0,\ghost}$, let $R_j$ be the subalgebra generated by $x_1,x_2$ and $z_k$, $k\in\I_{j,\ghost}$.
By \cite[Proposition 4.16]{AAH}, each $R_{j}$ is a free $R_{j+1}$-module with basis $(z_{j}^n)_{n\in\N_0}$.
Let $\sigma_j$ the algebra automorphism given by the action of $g_1^jg_2$, and $D_j:= [z_{j}, -]$. Then $D_j$ is a $(\id, \sigma_j)$-derivation, and by \cite[Lemmas 4.14, 4.15]{AAH}, 
$D_j$ restricts to a $(\id, \sigma_j)$-derivation of $R_{j+1}$. Thus $R_{j} \simeq R_{j+1}[z_{j};\sigma_j, D_j]$, for all $j\in \I_{0,\ghost}$.
Now $\toba(\lstr( 1, \ghost))$ is strongly noetherian
by \cite[Proposition 4.10]{ASZ};  AS-regular by \cite[Proposition 2]{AST} and Cohen-Macaulay by \cite[Lemma 5.3]{ZZ}.
\epf

\subsubsection{Organization of the Section}\label{subsubsec:generalities-organization}
In the rest of the Section we present the defining relations of some specific $\toba(Z_U)$. 
We state now some terminology for general use.
First, in \S  \ref{subsec:dim2} we deal with the case when $\dim W = 2$ and the braiding is of Cartan type.
The same method works routinely for other arithmetic root systems of rank 2.

Assume now that $\theta > 2$. Given $\Jb \subset \I_{\theta}$, we set
\begin{align*}
&\begin{aligned}
V_{\Jb} &= \oplus_{j\in \Jb} \ku x_j, & U_{\Jb} &= \oplus_{h\notin \Jb}\ku x_h,& &\text{so that}& W &= V_{\Jb}  \oplus U_{\Jb},
\end{aligned}\\
&Z_{U_{\Jb}}  = \ad_c\toba (V_{\Jb}) (U_{\Jb}) = \oplus_{h\notin \Jb} \ad \toba(V_{\Jb}) x_{h}.
\end{align*}
Observe that $Z_{V_{\Jb}} = Z_{U_{\Jb'}}$, where $\Jb' = \I_{\theta} - \Jb$ (thus we do not need to consider $Z_{V_{\Jb}}$)
and that $Z_{U_{\Jb}} \simeq Z_{U_{\vartheta(\Jb)}}$ as braided vector spaces whenever $\vartheta$ is an automorphism of the Dynkin diagram \cite[\S 2.1]{AA-diag-survey}.
We are mostly interested in the cases when $\Jb$ consists of an extremal vertex or its complement.
Indeed $Z_{U_{\Jb}}$ would turn out to be simple exactly when $\card \Jb = 1$, and $\toba (V_{\Jb})$
would give rise to $U_{q}(\g)$ with $\g$ simple exactly when the Dynkin diagram of $V_{\Jb}$ is connected.

\subsection{Dimension $W = 2$}\label{subsec:dim2}
Here $V = \ku x_1$ and $U = \ku x_2$ have dimension 1.
Then $\dim Z_U = n_{12} + 1$ and $\dim Z_V = n_{21} + 1$; more precisely,
\begin{align}
\notag Z_U  &= \ad_c\toba (V) (U) = \oplus_{n\in \I_{0, n_{12}}} \ku u_n,& u_n &:=\ad_c(x_1)^n(x_2);
\\ \label{eq:un}
u_{n+1} &= x_1 u_n - q_{11}^n q_{12} u_n x_1.
\end{align}
Notice that $u_n = x_{n\alpha_1 + \alpha_2}$ belongs to the PBW-basis alluded above.
Let
\begin{align*}
\acom{n} &= \prod_{j \in \I_{0,n-1}} (1 - q_{11}^j \widetilde{q}_{12}),& &\text{so that}
&  n \leq n_{12} &\implies  \acom{n} \neq 0.
\end{align*}

\begin{lema}
The braiding of $Z_U$ is given by
\begin{align}\label{eq:ZU-braid}
c(u_n \ot u_m) &= \sum_{j\in \I_{0,n}: m + j \leq n_{12}} q_{11}^{m(n-j)} q_{12}^{n-j}q_{21}^m q_{22}
\binom{n}{j}_{q_{11}} \frac{\acom{n}}{\acom{n-j}}  u_{m + j} \ot u_{n-j},
\end{align}
$ n,m \in \I_{0, n_{12}}$.  

In particular, $Z_U$ is of diagonal type if and only if $n_{12} = 0$, hence also $n_{21} = 0$.
\end{lema}
\pf
Recall the realization in \S  \ref{subsec:generalities}.
The coaction $\delta: Z_U \to \toba (V)\# \ku \Gamma \ot Z_U$ is given by \eqref{eq:coaction-R}. Hence
$\delta(u_0) = \alpha_2 \ot u_0$. 
We claim that for every $n \in \I_{0, n_{12}}$,
\begin{align}\label{eq:delta-un}
\delta(u_n) &= \sum_{i\in \I_{0,n}} \binom{n}{i}_{q_{11}} \frac{\acom{n}}{\acom{i}}  x_1^{n-i} g_1^i g_2 \ot u_i. 
\end{align}
The proof of  \eqref{eq:delta-un} goes by induction on $n$. Indeed, by \eqref{eq:un}, we have
\begin{align*}
\delta(u_{n+1}) = (x_1 \ot 1 + g_1 \ot x_1) \delta(u_n) - q_{11}^n q_{12} \delta(u_n) (x_1 \ot 1 + g_1 \ot x_1).
\end{align*}
Then \eqref{eq:ZU-braid} follows. 
We next claim that 
\begin{align}\label{eq:c-unum}
\exists q \in \ku^{\times}, x \in Z_U:\ c(x\ot x) &= q x\ot x \iff x \in \ku u_0 \text{ or } x \in \ku u_{n_{12}}.
\end{align}

Assume that there exist such $q$ and $x$ and write
\begin{align*}
x &= \sum_{i\in \I_{n,p}} c_i u_i,&   n \le p, \, c_n, c_p &\neq 0.
\end{align*}
Let $Q^{i,j}_t = q_{11}^{j(i-t)} q_{12}^{i-t}q_{21}^j q_{22}
\binom{i}{t}_{q_{11}} \frac{\acom{i}}{\acom{i-t}}$. 
Then
\begin{align*}
\sum_{\substack{ i,j\in \I_{n,p}, \\ t\in \I_{0,i}: j + t \leq n_{12}}} c_i c_j Q^{i,j}_t  u_{j + t} \ot u_{i-t}
 &=  \sum_{r,s\in \I_{n,p}} q c_r c_s u_r \ot u_s;
 \\
\implies
\sum_{\substack{ i,j\in \I_{n,p}, \\ 0 \le t \leq \min\{i, n_{12} - j\}\\ j+t = r, i-t =s }} c_i c_j Q^{i,j}_t  u_{j + t} \ot u_{i-t} &= q c_r c_s u_r \ot u_s
\quad \forall r,s\in \I_{n,p}.
\end{align*}
Take $(r, s)=(n,p)$. Then $j=n$, $i=p$ and necessarily $t=0$; the last could only happen if $\min\{p, n_{12} - n\} = 0$.
That is, either $p=0=n$, or else $n = n_{12} = p$. This shows \eqref{eq:c-unum}. 

Thus, if $Z_U$ is diagonal type, then $n_{12} \leq 1$. If $n_{12} = 1$, then  
\begin{align*}
c(u_0 \otimes u_1) &= q_{21} q_{ 22} u_1 \otimes u_0, \\ c(u_1 \otimes u_0) &= q_{12} q_{22} u_0 \otimes u_1 + q_{22} (1 - \widetilde{q}_{12}) u_1 \otimes u_0
\end{align*}
by a similar computation. Here $\widetilde{q}_{12} \neq 1$ since $n_{12} \neq 0$; hence $Z_U$ is not of diagonal type.
We conclude that $n_{12} = 0$.
\epf

\subsubsection{Cartan type $A_2$}\label{exa:A2} Here $q \in \ku^{\times} - \{1\}$. The Dynkin diagram is $\Dynkintwoxy{q}{1}{q^{-1}}{q}{2}$ and \begin{align*}
\varDelta_+ = \{\beta_1 = \alpha_1, \beta_2 = \alpha_1 + \alpha_2, \beta_3 = \alpha_2 \}.
\end{align*}
Thus $Z_U = \ku u_0 \oplus \ku u_1$; $u_0 = x_{\beta_{3}}$, $u_1 = x_{\beta_{2}}$. By Remark \ref{rem:consecutive-roots}, $[u_1, u_0]_c = 0$. Hence 
\begin{align*}
\toba(Z_U) = \ku \langle u_0, u_1\vert u_1u_0 - q_{12} q\, u_0 u_1, u_0^N, u_1^N \rangle,
\end{align*}
where $N = \ord q$. We omit the last relations when $N = \infty$, in which case:

\begin{prop}
If $W$ is of type $A_2$ and $q \notin \G_{\infty}$, then $\toba(Z_U)$ is a quantum plane, with $\GK = 2$.
\end{prop}

It is well-known that in this case, $\toba(Z_U)$ is AS-regular, see \cite[(0.2)]{ArS}.

\subsubsection{Cartan type $B_2$}\label{exa:B2} The Dynkin diagram is $\Dynkintwoxy{q}{1}{q^{-2}}{q^2}{2}$, 
$\ord q > 2$, and 
\begin{align*}
\varDelta_+ = \{\beta_1 = \alpha_1, \beta_2 = 2\alpha_1 + \alpha_2, \beta_3= \alpha_1 + \alpha_2, \beta_4 = \alpha_2 \}.
\end{align*}
Hence $Z_U = \ku u_0 \oplus \ku u_1 \oplus \ku u_2$; $u_i = x_{\beta_{4-i}}$, $i\in \I_{0,2}$. 
By Remark \ref{rem:consecutive-roots}, 
\begin{align}\label{eq:consecutive-roots-B2}
u_1u_0 &= q^2q_{12} \, u_0u_1, &
u_2u_1 &= q^2q_{12} \, u_1u_2.
\end{align}
We check that
\begin{align}\label{eq:other-relations-B2}
u_2u_0 &= q^2q_{12}^2 \, u_0u_2 + qq_{12}(q-1) \, u_1^2.
\end{align}
Thus $\toba(Z_U)$ is presented by generators $u_i$, $i\in\I_{0,2}$, and relations \eqref{eq:consecutive-roots-B2}, \eqref{eq:other-relations-B2}, 
\begin{align*}
u_0^M&=0, & u_1^N&=0, & u_2^M&=0,
\end{align*}
where $N = \ord q$, $M = \ord q^2$. Particularly, when $N = \infty$, we have:

\begin{prop}
If $W$ is of type $B_2$ and $q \notin \G_{\infty}$, then $\toba(Z_U)$ is a quadratic algebra with relations \eqref{eq:consecutive-roots-B2}, 
\eqref{eq:other-relations-B2}, and $\GK = 3$.
\end{prop}

In this case, it is known that  $\toba(Z_U)$ is AS-regular, see \cite[Remark 2.9]{NS}.

\subsubsection{Cartan type $C_2$}\label{exa:C2} The Dynkin diagram is $\Dynkintwoxy{q^2}{1}{q^{-2}}{q}{2}$, $\ord q > 2$, and \begin{align*}
\varDelta_+ = \{\beta_1 = \alpha_1, \beta_2 = \alpha_1 + \alpha_2, \beta_3= \alpha_1 + 2\alpha_2, \beta_4 = \alpha_2 \}.
\end{align*}
Hence $Z_U = \ku u_0 \oplus \ku u_1$, $u_0 = x_{\beta_{4}}$, $u_1 = x_{\beta_{2}}$. Set $u_{10}=u_1u_0-q_{12}q u_0u_1 =x_{\beta_3}$. By Remark \ref{rem:consecutive-roots}, 
\begin{align}\label{eq:consecutive-roots-C2}
u_1u_{10} &= q^2q_{12} \, u_{10}u_1, &
u_{10}u_0 &= q^2q_{12} \, u_0u_{10}.
\end{align}
Thus $\toba(Z_U)$ is presented by generators $u_0$, $u_1$, and relations \eqref{eq:consecutive-roots-C2}, 
\begin{align*}
u_0^N&=0, & u_{10}^M&=0, & u_1^N&=0,
\end{align*}
where $N = \ord q$, $M = \ord q^2$. Clearly \eqref{eq:consecutive-roots-C2} are quantum Serre relations.
\begin{obs}
If $N < \infty$, then $\toba(Z_U)$ 
is isomorphic as an algebra to a Nichols algebra of diagonal type only when $M = N$.
\end{obs}

Particularly, when $N = \infty$, we have:

\begin{prop}
If $W$ is of type $C_2$ and $q \notin \G_{\infty}$, then $\toba(Z_U)$ is a cubic algebra with relations \eqref{eq:consecutive-roots-C2}, and $\GK = 3$.
\end{prop}

In this case, it is known that  $\toba(Z_U)$ is AS-regular, see \cite[(8.5)]{ArS}.

\subsubsection{Cartan type $G_2$}\label{exa:G2} The Dynkin diagram is $\Dynkintwoxy{q}{1}{q^{-3}}{q^3}{2}$, $\ord q > 3$, and 
\begin{align*}
\varDelta_+ &= \{\beta_1 = \alpha_1, \beta_2 = 3\alpha_1 + \alpha_2, \beta_3= 2\alpha_1 + \alpha_2, \\ & \quad \beta_4 = 3\alpha_1+2\alpha_2, \beta_5= \alpha_1 + \alpha_2, \beta_6 = \alpha_2 \}.
\end{align*}
Hence $Z_U = \ku u_0 \oplus \ku u_1 \oplus \ku u_2 \oplus \ku u_3$, $u_0 = x_{\beta_{6}}$, $u_1 = x_{\beta_{5}}$, $u_2 = x_{\beta_{3}}$, $u_3 = x_{\beta_{2}}$. 
Set $u_{21}=u_2u_1 -q^2q_{12} u_1u_2 =x_{\beta_4}$.
By Remark \ref{rem:consecutive-roots}, 
\begin{align}\label{eq:consecutive-roots-G2}
u_1u_0 &= q^3q_{12} \, u_0u_1, &
u_3u_2 &= q^3q_{12} \, u_2u_3.
\end{align}
We check that
\begin{align}\label{eq:other-relations-G2}
&\begin{aligned}
u_3u_1 &= q^3q_{12}^2 \, u_1u_3 + \frac{qq_{12}(q^3-1)}{q+1} \, u_2^2, \\
u_3u_0 &= q^3q_{12}^3 \, u_0u_3 + q_{12}(q^3-q^2-q) \, u_{21} + qq_{12}(1-q)(3)_q^{!} \, u_{1}u_2, \\
u_2u_0 &= q^3q_{12}^2 \, u_0u_2 + qq_{12}(q^2-1) \, u_1^2, 
\end{aligned}
\\\label{eq:other-relations-G2bis}
&\begin{aligned}
u_2u_{21} &= q^3q_{12} \, u_{21}u_2, \\ u_{21} u_1 &= q^3q_{12} \, u_1u_{21}.
\end{aligned}
\end{align}
Thus $\toba(Z_U)$ is generated by $u_i$, $i\in\I_{0,3}$, with relations \eqref{eq:consecutive-roots-G2}, \eqref{eq:other-relations-G2}, \eqref{eq:other-relations-G2bis},
\begin{align*}
u_0^M&=0, & u_1^N&=0, & u_{21}^M&=0, & u_2^N&=0, & u_3^M&=0,
\end{align*}
where $N = \ord q$, $M = \ord q^3$. 
Particularly, when $N = \infty$, we have:

\begin{prop}
If $W$ is of type $G_2$ and $q \notin \G_{\infty}$, then $\toba(Z_U)$ has quadratic relations \eqref{eq:consecutive-roots-G2}, \eqref{eq:other-relations-G2}
and cubic relations \eqref{eq:other-relations-G2bis}; also, $\GK = 5$.
\end{prop}

\subsubsection{Cartan type $G_2$ with the inverse numeration}\label{exa:G2-2} The Dynkin diagram is $\Dynkintwoxy{q^3}{1}{q^{-3}}{q}{2}$, 
$\ord q > 3$, and
\begin{align*}
\varDelta_+ &= \{\beta_1 = \alpha_1, \beta_2 = \alpha_1 + \alpha_2, \beta_3= 2\alpha_1 + 3\alpha_2, 
\\ & \quad \beta_4 = \alpha_1+2\alpha_2, \beta_5= \alpha_1 + 3\alpha_2, \beta_6 = \alpha_2 \}.
\end{align*}
Hence $Z_U = \ku u_0 \oplus \ku u_1$, $u_0 = x_{\beta_{6}}$, $u_1 = x_{\beta_{2}}$. Set $u_{10}=u_1u_0 -qq_{12} u_0u_1 =x_{\beta_4}$,
$u_{100}=u_{10}u_0 -q^2q_{12} u_0u_{10} = x_{\beta_5}$, $u_{110}=u_1u_{10} -q^2q_{12} u_{10}u_1 =x_{\beta_3}$.
By Remark \ref{rem:consecutive-roots}, 
\begin{align}\label{eq:consecutive-roots-G2-2}
u_1u_{110} &= q^3q_{12} \, u_{110}u_1, &
u_{100}u_0 &= q^3q_{12} \, u_0u_{100}.
\end{align}
We check that
\begin{align}\label{eq:other-relations-G2-2}
u_1u_{100} &= q^3q_{12}^2 \, u_{100}u_1 + \frac{q_{12}(3)_q(q-1)^2}{q+1} \, u_{10}^2.
\end{align}
Thus $\toba(Z_U)$ is presented by generators $u_0$, $u_1$, and relations \eqref{eq:consecutive-roots-G2-2}, \eqref{eq:other-relations-G2-2}, 
\begin{align*}
u_0^N&=0, & u_{110}^M&=0, & u_{10}^N&=0, & u_{100}^M&=0, & u_1^N&=0,
\end{align*}
where $N = \ord q$, $M = \ord q^3$. When $N = \infty$, we have:

\begin{prop}
If $W$ is of type $G_2$ and $q \notin \G_{\infty}$, then $\toba(Z_U)$ is defined by the
quartic relations \eqref{eq:consecutive-roots-G2-2}, \eqref{eq:other-relations-G2-2}; also, $\GK = 5$.
\end{prop}

In this case, $\toba(Z_U)$ is AS-regular; it is the \emph{Algebra F} in \cite[\S 4]{WW}.

\subsection{Cartan type $A_\theta$, $\theta > 2$}\label{subsec:An} Here $q \in \ku^{\times} - \{1\}$.
The Dynkin diagram is 
\begin{align*}
\xymatrix{  \overset{q}{\underset{1}{\circ}}\ar  @{-}[r]^{q^{-1}}  &
	\overset{q}{\underset{2}{\circ}} \ar@{.}[r] & \overset{q}{\underset{\theta - 1}{\circ}} \ar  @{-}[r]^{q^{-1}}  
	& \overset{q}{\underset{\theta}{\circ}}}.
\end{align*}
With the notation \eqref{eq:alfaij} and ordered lexicographically, the positive roots are
\begin{align} \label{eq:Atheta-roots}
\varDelta_+ = \{\alpha_{k\, j}\,|\, k \leq j \in \I\}.
\end{align}
The generators of the PBW-basis are $x_{\alpha_{kk}} = x_k$, $x_{\alpha_{kl}} = x_{(kl)}$, $k < l$, cf. \eqref{eq:roots-Atheta}.
Let $N = \ord q$.
Assume that $q \neq -1$. The defining relations are
\begin{align}\label{eq:rels-serre-conmutan}
x_{ij} &= 0,& i &< j - 1, 
\\  \label{eq:rels-serre-2}
x_{iij} &= 0,& \vert j - i\vert &= 1, 
\\ \label{eq:rels-prv}
x_{(kl)}^N &=0,& k &\leq l,
\end{align}
where $i,j,k,l \in \I$. When $q = -1$, the defining relations are  \eqref{eq:rels-serre-conmutan}, \eqref{eq:rels-prv} (with $N =2$), and
\begin{align} \label{eq:rels-serre-2-menos1}
[x_{(i-1, i+1)}, x_i]_c &= 0,
\end{align}
$i\in \I_{2,\theta - 1}$. As usual we omit \eqref{eq:rels-prv} when $q\notin \gi$.

\subsubsection{$\Jb = \{1\}$}\label{exa:Atheta-1}
The algebra  $\toba(Z_U)$
is  presented by $x_{12}, x_{2}, x_{3},  \dots,  x_{\theta}$ with defining relations 
\eqref{eq:rels-serre-conmutan}, $i, j \in \I_{2, \theta}$; 
\eqref{eq:rels-serre-2} if $q\neq -1$, or \eqref{eq:rels-serre-2-menos1} if $q= -1$, $i, j \in \I_{2, \theta}$;
 \eqref{eq:rels-prv}, $k,l \in \I_{2, \theta}$; 
and the new relations
\begin{align}\label{eq:rels-Zu-conmutan}
[x_{12}, x_i]_c &= 0,& i &\in \{2\} \cup \I_{4, \theta}, 
\\\label{eq:rels-Zu-serre}
[x_{12}, [x_{12}, x_3]_c]_c &= 0,&  x_{3312} &=0, 
\\\label{eq:rels-Zu-serre-gen}
[[x_{12}, x_3]_c, x_2]_c &= 0,& & 
\\ \label{eq:rels-prv12-Zu}
x_{12}^N &=0,&  &
\\ \label{eq:rels-prv-Zu}
([x_{12}, x_{(3 l)}]_c)^N &=0,& l &\in \I_{3, \theta}.
\end{align}

We omit the last relations when $N = \infty$, in which case $\toba(Z_U)$ is a domain and $\GK \toba(Z_U) =  \binom{\theta + 1}{2} - 1$.

\subsubsection{$\card\Jb = 1$}\label{exa:A3-2}
The case $\Jb = \{\theta\}$ reduces to the previous one.
To exemplify, we discuss only the case $\theta = 3$, $\Jb = \{2\}$. 
The algebra  $\toba(Z_U)$
is  presented by $x_{1}, x_{21}, x_{23},  x_{3}$ with defining relations 
\begin{align}\label{eq:rels-Zu-2-conmutan}
&\begin{aligned}
&[x_{1}, x_{3}]_c = 0, & [x_{1}, x_{21}]_c &= 0, & [x_{23}, x_{3}]_c &= 0, & [x_{21}, x_{23}]_c&=0,
\end{aligned}
\\ 
\label{eq:rels-Zu-2-serre}
&\begin{aligned}
&[x_{21}, x_{3}]_c = q_{13}(1-q)x_{23}x_1 -q_{21} [x_{1}, x_{23}]_c,   
\end{aligned}
\\ \label{eq:rels-2-prv12-Zu}
&\begin{aligned}
x_{12}^N &=0,&  x_{21}^N&=0,& x_{23}^N&=0,& x_{3}^N&=0, & [x_{1}, x_{23}]_c^N&=0.
\end{aligned}
\end{align}

\subsubsection{$\Jb = \I_{\theta - 1}$}\label{exa:A2-theta-1}

The algebra  $\toba(Z_U)$
is  presented by $x_{(1\theta)}, x_{(2\theta)},  \dots,  x_{\theta}$ with defining relations 
\begin{align}\label{eq:rels-Zu-qls}
[x_{(i\theta)}, x_{(j\theta)}]_c &= 0,& i &< j\in  \I_{\theta}, 
\\ \label{eq:rels-prv-qls-Zu}
x_{(i\theta)}^N &=0,&  i &\in  \I_{\theta}.
\end{align}

We omit the last relations when $N = \infty$, in which case $\toba(Z_U)$ is a quantum linear space and $\GK \toba(Z_U) =  \theta$.

\subsection{Cartan type $B_\theta$, $\theta > 2$}\label{subsec:Bn} 
Here $q \in \ku^{\times} - \{\pm 1\}$.
The Dynkin diagram is 
\begin{align*}
\xymatrix{ \overset{\,\,q^2}{\underset{\ }{\circ}}\ar  @{-}[r]^{q^{-2}}  &
	\overset{\,\,q^2}{\underset{\ }{\circ}}\ar  @{-}[r]^{q^{-2}} &
	\overset{\,\,q^2}{\underset{\ }{\circ}}\ar@{.}[r] & \overset{\,\,q^2}{\underset{\
		}{\circ}} \ar  @{-}[r]^{q^{-2}}  & \overset{q}{\underset{\ }{\circ}}}.
\end{align*}
With the notation \eqref{eq:alfaij} and ordered lexicographically, the positive roots are
\begin{align} \label{eq:root-system-B}
	\varDelta^\bq_+&=\{\alpha_{ij}\,|\, i\leq j \in\I \} \cup 
	\{\alpha_{i\theta} + \alpha_{j\theta}\,|\, i< j\in\I \}.
\end{align}
The generators of the PBW-basis are $x_{\alpha_{kk}} = x_k$, $x_{\alpha_{kl}} = x_{(kl)}$, $k < l \in \I_{\theta}$,
$x_{\alpha_{i\theta} + \alpha_{\theta}} = [x_{\alpha_{i\theta}}, x_\theta]_c$, $i  \in \I_{\theta - 1}$,
$x_{\alpha_{i\theta} + \alpha_{j\theta}} = [x_{\alpha_{i\theta} + \alpha_{(j+1) \theta}}, x_j]_c$, 
$i <  j \in \I_{\theta - 1}$. 

For simplicity, we assume that either $N > 4$ is odd or else is $\infty$. 
The defining relations are
\begin{align}\label{eq:rels-type-B-N>4-odd}
x_{ij} &= 0,& i &< j - 1; & x_{iii\pm1}&=0, & i &< \theta;
\\ \label{eq:rels-type-B-N>4-odd-2}
x_{\theta\theta\theta\theta-1}&=0; 
\\ \label{eq:prv-B}
x_{\alpha}^{N} &=0,& &	\alpha\in\varDelta_{+}. 
\end{align} 
where $i,j \in \I$. See \cite[\S 4.2]{AA-diag-survey}
for the relations in other cases.
As usual we omit \eqref{eq:prv-B} when $q\notin \gi$.

\subsubsection{$\Jb = \{1\}$}\label{exa:Btheta-1}
The algebra  $\toba(Z_U)$
is  presented by $x_{12}, x_{2}, x_{3},  \dots,  x_{\theta}$ with defining relations 
\eqref{eq:rels-type-B-N>4-odd}, $i, j \in \I_{2, \theta}$; \eqref{eq:rels-type-B-N>4-odd-2};
\eqref{eq:prv-B}, $\supp \alpha \subset \I_{2, \theta}$; 
\eqref{eq:rels-Zu-conmutan}; \eqref{eq:rels-Zu-serre}; \eqref{eq:rels-Zu-serre-gen}; \eqref{eq:rels-prv12-Zu}; \eqref{eq:rels-prv-Zu}; and the new relations
\begin{align}\label{eq:rels-prv-Zu-B}
([\cdots [[x_{12}, x_{(3\theta)}]_c,x_{\theta} ]_c, \cdots, \ x_l]_c  )^N &=0,& l &\in \I_{2, \theta}.
\end{align}

\subsubsection{$\Jb = \{\theta\}$}\label{exa:Btheta-theta}
The algebra  $\toba(Z_U)$
is  presented by $x_{1}, x_{2}, \dots,  x_{\theta-1}$, $u_1=x_{\theta \, \theta-1}$, $u_2=x_{\theta\theta \, \theta-1}$, with defining relations 
\eqref{eq:rels-type-B-N>4-odd}, $i, j \in \I_{\theta}$; 
\eqref{eq:prv-B}, $\alpha= \alpha_{kl}$, $k<l\in\I_{\theta-1}$; 
and the new relations
\begin{align} \label{eq:rels-ZU-B-theta-conmutan}
&\begin{aligned}
& [x_i,u_1]_c =0, & [x_i,u_2]_c &=0, & i &\in\I_{\theta-3};
\end{aligned}
\\ \label{eq:rels-ZU-B-theta-u1}
&\begin{aligned}
& [[x_{\theta-2},u_1]_c ,u_1]_c = -q_{\theta \,\theta-1}[x_{\theta-2 \, \theta-1},u_2]_c;
\end{aligned}
\\ \label{eq:rels-ZU-B-theta-u1u2}
&\begin{aligned}
&[[x_{\theta-2},u_2]_c ,u_1]_c=0, & [[x_{\theta-2},u_2]_c ,u_2]_c &=0;
\end{aligned}
\\ \label{eq:rels-ZU-B-theta-PRV1}
&\begin{aligned}
&[x_{(i \, \theta-2)},u_1]_c^N =0, & [x_{(i \, \theta-2)},u_2]_c^N &=0, & i &\in\I_{\theta-2};
\end{aligned}
\\ \label{eq:rels-ZU-B-theta-PRV2}
&\begin{aligned}
&[\dots[[x_{(i \, \theta-2)},u_2]_c, x_{\theta-1}]_c, \dots x_j]_c^N =0, & i<j &\in\I_{\theta-2};
\end{aligned}
\\ \label{eq:rels-ZU-B-theta-PRV3}
&\begin{aligned}
&u_1^N=0; & u_2^N &=0.
\end{aligned}
\end{align}

\subsubsection{$\Jb = \I_{\theta -1}$}\label{exa:Btheta-1-thetamenos1}

Here, $Z_U$ is spanned by $z_i:=x_{(i\theta)}$, $i\in\I$.

The algebra $\toba(Z_U)$ is  presented by $z_i$, $i\in\I$, with defining relations
\begin{align} \label{eq:rels-ZU-Btheta-1-thetamenos1-qsr}
&\begin{aligned}
& [z_i,[z_i,z_j]_c]_c =0, & &[[z_i,z_j]_c,z_j]_c=0, & &i<j\in\I;
\\ 
& [[z_i,z_j]_c, z_k]_c =0, & &[[z_i,z_k]_c,z_j]_c=0, & & i<j<k \in\I;
\end{aligned}
\\ \label{eq:rels-ZU-Btheta-1-thetamenos1-PRV1}
&\begin{aligned}
&z_i^N =0, \ i\in\I; &
&[z_i, z_j]_c^N =0, \ i<j\in\I.
\end{aligned}
\end{align}

\subsubsection{$\Jb = \I_{2, \theta}$}\label{exa:Btheta-2-theta}

Here, $Z_U$ is spanned by 
\begin{align*}
w_i&:=x_{i \, i-1 \dots 1}, \ i\in\I, && \text{and} & \widetilde{w}_j &:=x_{j \dots \theta \, \theta \dots 1}, \ j\in\I_{2,\theta}.
\end{align*}

The algebra $\toba(Z_U)$ is  presented by $w_i$, $\widetilde{w}_i$, $i\in\I$, with defining relations
\begin{align} \label{eq:rels-ZU-Btheta-2-theta-qsr}
&\begin{aligned}
& [w_i,w_j]_c=0, \ i<j\in\I; &  &[\widetilde{w}_i, \widetilde{w}_j]_c = 0, \ i<j\in\I_{2,\theta};
\\ 
& [\widetilde{w}_j,w_{j-1}]_c = q_{j\, j-1}q (q-1) w_{\theta}^2, & & j\in\I_{2,\theta}; & 
\\ 
&[\widetilde{w}_j,w_i]_c=0,& &i\in\I, j\in\I_{2,\theta}-\{i+1\};
\end{aligned}
\\ \label{eq:rels-ZU-Btheta-2-theta-PRV1}
&\begin{aligned}
&w_i^N =0, \ i\in\I; &
& \widetilde{w}_j^N =0, \ j\in\I_{2,\theta}.
\end{aligned}
\end{align}

\subsection{Cartan type $C_\theta$, $\theta > 2$}\label{subsec:Cn} 
Here $q \in \ku^{\times} - \{\pm 1\}$.
The Dynkin diagram is 
\begin{align*}
\xymatrix{ \overset{\,\,q}{\underset{\ }{\circ}}\ar  @{-}[r]^{q^{-1}}  &
	\overset{\,\,q}{\underset{\ }{\circ}}\ar  @{-}[r]^{q^{-1}} &
	\overset{\,\,q}{\underset{\ }{\circ}}\ar@{.}[r] & \overset{\,\,q}{\underset{\
		}{\circ}} \ar  @{-}[r]^{q^{-2}}  & \overset{q^2}{\underset{\ }{\circ}}}.
\end{align*}
The set of positive roots is
\begin{align}\label{eq:root-system-C}
\varDelta^+&=\{\alpha_{i\, j}\,|\, i\leq j\in\I \}\cup \{\alpha_{i\, \theta}+\alpha_{j\, \theta-1}\,|\, i\leq j\in\I_{\theta-1} \}.
\end{align}
The generators of the PBW-basis are $x_{\alpha_{kk}} = x_k$, 
$x_{\alpha_{ij}}=x_{(ij)}$, $i <  j \in \I$, 
$x_{\alpha_{i\theta} + \alpha_{i\theta-1}}= [x_{(i\theta)},x_{(i\theta-1)}]_c$, $x_{\alpha_{i\theta} + \alpha_{\theta-1}}=[x_{(i\theta)}, x_{\theta-1}]_c$,  $i \in \I_{\theta-1}$, 
$x_{\alpha_{i\theta} + \alpha_{j\theta-1}} = [x_{\alpha_{i\theta} + \alpha_{j+1\theta-1}}, x_j]_c$, $i <  j \in \I_{\theta-2}$.
For simplicity, we assume that either $N > 4$ is odd or else is $\infty$.  The defining relations are
\begin{align}\label{eq:rels-type-C-N>3-qsr}
&x_{ij} = 0, \quad i < j - 1; &  x_{iij}&=0, \quad j=i\pm 1, (i,j)\neq(\theta-1,\theta);
\\ \label{eq:rels-type-C-N>3-qsr-bis}
&x_{iii\theta}=0, \ i = \theta - 1; &&
\\ \label{eq:rels-type-C-N>3}
&x_{\alpha}^{N} =0, &&\alpha \in\varDelta_{+}.
\end{align}
See \cite[\S 4.2]{AA-diag-survey} for the relations in other cases. As usual we omit \eqref{eq:rels-type-C-N>3} when $q\notin \gi$.

\subsubsection{$\Jb = \I_{\theta -1}$}\label{exa:Ctheta-1-thetamenos1}

Here, $Z_U$ is spanned by 
\begin{align*}
z_i&:=x_{(i\theta)}, \ i\in\I, & &\text{and} & y_{ij}:= [x_{(i\theta-1)},x_{(j\theta)}]_c, \ i\le j \in\I_{\theta-1}.
\end{align*}
The algebra $\toba(Z_U)$ is  presented by generators $z_i$, $i\in\I$, and $y_{ij}$, $i<j \in\I_{\theta-1}$, with defining relations
\begin{align} \label{eq:rels-ZU-Ctheta-1-thetamenos1-qsr}
&\begin{aligned}
& [z_i,z_j]_c=0, & &i<j\in\I;
\\ 
& [y_{ii},z_i]_c=0, \quad  [y_{ii},z_{\theta}]_c=\bq_{\alpha_{(i\theta)} \alpha_{\theta}}(1-q^{-1}) z_{i}^2, & &i\in\I_{\theta-1};
\\
& [y_{ij},z_{\theta}]_c=\bq_{\alpha_{(i\theta-1)} \alpha_{(j\theta)}}(q^{2}-1) z_{j}z_{i}, & &i<j\in\I_{\theta-1};
\\
& [y_{ij},z_{k}]_c=\bq_{\alpha_{j\theta} \alpha_{k \theta}} (q-1) z_{j}y_{ik}, & &i\le j<k\in\I_{\theta-1};
\\
& [y_{ik},z_{j}]_c=0, \qquad [z_{i},y_{jk}]_c=0, & &i\le j \le k\in\I_{\theta-1};
\\
& [y_{ij},y_{kl}]_c=0, & &i,j,k,l \in\I_{\theta-1};
\end{aligned}
\\ \label{eq:rels-ZU-Ctheta-1-thetamenos1-PRV1}
&\begin{aligned}
&z_i^N =0, \ i\in\I; &
&y_{ij}^N =0, \ i\le j\in\I_{\theta-1}.
\end{aligned}
\end{align}

\subsubsection{$\Jb = \I_{2, \theta}$}\label{exa:Ctheta-2-theta}

Here, $Z_U$ is spanned by 
\begin{align*}
w_i&:=x_{i \, i-1 \dots 1}, \ i\in\I, && \text{and} & \widetilde{w}_j &:=[x_{\theta \dots j}, w_{\theta-1}]_c, \ j\in\I_{\theta-1}.
\end{align*}

The algebra $\toba(Z_U)$ is  presented by generators $w_i$, $\widetilde{w}_i$, $i\in\I$, with defining relations
\begin{align} \label{eq:rels-ZU-Ctheta-2-theta-qsr}
&\begin{aligned}
& [w_i,w_j]_c=0, \quad i<j\in\I, \ (i,j)\neq (\theta-1,\theta);
\\
& [w_{\theta-1},[w_{\theta-1},w_{\theta}]_c]_c= [w_{\theta},[w_{\theta},w_{\theta-1}]_c]_c = 0
\\  
&[\widetilde{w}_i, \widetilde{w}_j]_c = 0, \quad i<j\in\I_{\theta-1};
\\ 
& [\widetilde{w}_j,w_{j-1}]_c = \bq_{\alpha_{(j\theta)} \alpha_{(1 \ \theta-1)}}(q-1) w_{\theta}w_{\theta-1}, \quad j\in\I_{\theta-1}; 
\\ 
&[\widetilde{w}_j,w_i]_c=0, \quad i\in\I, j\in\I_{\theta-1}-\{i+1\};
\end{aligned}
\\ \label{eq:rels-ZU-Ctheta-2-theta-PRV1}
&\begin{aligned}
&w_i^N =0, \ i\in\I; &
& \widetilde{w}_j^N =0, \ j\in\I_{\theta-1}; & & [w_{\theta},w_{\theta-1}]_c^N=0.
\end{aligned}
\end{align}

\subsection{Cartan type $D_\theta$, $\theta > 3$}\label{subsec:Dn} 
Here $q \in \ku^{\times} - \{1\}$.
The Dynkin diagram is 
\begin{align}\label{eq:dynkin-type-D}
\xymatrix{ & & & &  \overset{q}{\circ} &\\
	\overset{q}{\underset{\ }{\circ}}\ar  @{-}[r]^{q^{-1}}  & \overset{q}{\underset{\
		}{\circ}}\ar  @{-}[r]^{q^{-1}} &  \overset{q}{\underset{\ }{\circ}}\ar@{.}[r] &
	\overset{q}{\underset{\ }{\circ}} \ar  @{-}[r]^{q^{-1}}  & \overset{q}{\underset{\
		}{\circ}} \ar @<0.7ex> @{-}[u]_{q^{-1}}^{\qquad} \ar  @{-}[r]^{q^{-1}} &
	\overset{q}{\underset{\ }{\circ}}}
\end{align}
With the notation \eqref{eq:alfaij} and ordered lexicographically, the positive roots are
\begin{align}\label{eq:root-system-D}
\begin{aligned}
\varDelta^+&=\{\alpha_{i\, j}\,|\, i\leq j\in\I, \, (i,j)\neq (\theta-1,\theta) \}
\\ 
& \qquad \cup \{\alpha_{i\, \theta-2}+\alpha_{\theta}\,|\, i\in\I_{\theta-2} \} \cup \{\alpha_{i\, \theta}+\alpha_{j\, \theta-2}\,|\, i<j\in\I_{\theta-2} \}.
\end{aligned}\end{align}
The generators of the PBW-basis are 
\begin{align*}
x_{\alpha_{kk}} &= x_k,\quad k \in \I, &
x_{\alpha_{ij}} &= x_{(ij)},\quad i <j \in \I_{\theta-1},\\
x_{\alpha_{i\theta-2}+\alpha_{\theta}} &= [x_{(i\theta-2)}, x_{\theta}]_c,&
x_{\alpha_{i\theta}}&= [x_{\alpha_{i\theta-2}+\alpha_{\theta}}, x_{\theta-1}]_c, \quad i\in \I_{\theta-2},\\
x_{\alpha_{i\theta} + \alpha_{j\theta-2}} &= [x_{\alpha_{i\theta} + \alpha_{j+1\theta-2}}, x_j]_c,&
 i &<  j \in \I_{\theta-2}.
\end{align*}
For simplicity, we assume that either $N > 2$ or else is $\infty$.  The defining relations are
\begin{align}\label{eq:rels-type-D-N>2}
\begin{aligned}
x_{(\theta-1)\theta} &=0; & x_{ij} &= 0, \quad i < j - 1, (i,j)\neq (\theta-2,\theta); \\
x_{ii\theta}&=0, \ i= \theta-2; & x_{iij} &= 0, \quad \vert j - i\vert = 1, i,j\neq \theta;
\\
x_{ii(\theta-2)}&=0, \ i= \theta; & x_{\alpha}^{N} &=0,
\quad \alpha\in\varDelta_{+}.
\end{aligned}
\end{align}
See \cite[\S 4.2]{AA-diag-survey} for the relations in other cases. As usual we omit the last set of relations when $q\notin \gi$.

\subsubsection{$\Jb = \I_{\theta -1}$}\label{exa:Dtheta-1-thetamenos1}

Here, $Z_U$ is spanned by $x_{\theta}$,
\begin{align*}
z_i&:=[x_{(i\theta-2)},x_{\theta}], & \widetilde{z}_i &:= [x_{\theta-1},z_i]_c, & i & \in\I, \\ y_{ij} &:= [x_{(i\theta-1)},z_j]_c, & i < j & \in\I_{\theta-2}.
\end{align*}
The algebra $\toba(Z_U)$ is  presented by generators $x_{\theta}$, $z_i$, $\widetilde{z}_i$, $i\in\I$, and $y_{ij}$, $i<j \in\I_{\theta-1}$, with defining relations
\begin{align} \label{eq:rels-ZU-Dtheta-1-thetamenos1-qsr}
&\begin{aligned}
& [z_i,z_j]_c=0, \quad [\widetilde{z}_i,\widetilde{z}_j]_c=0, \quad [\widetilde{z}_i,z_j]_c=0, & &i<j\in\I_{\theta-2};
\\ 
&  [z_{i},\widetilde{z}_{j}]_c=\bq_{\alpha_{(ij-1)}, \alpha_{(j\theta-2)}+\alpha_\theta}(q-1)z_j\widetilde{z}_{i}, & &i<j\in\I_{\theta-2};
\\
& [z_i,x_{\theta}]_c=0, \quad [\widetilde{z}_i,x_{\theta}]_c=0, & & i\in\I_{\theta-2};
\\ 
& [y_{ij},z_{\theta}]_c=\bq_{\alpha_{(i\theta-1)}, \alpha_{(j\theta-2)}+\alpha_{\theta}}(q-1) z_{j}\widetilde{z}_{i}, & &i<j\in\I_{\theta-2};
\\
& [y_{ij},z_{k}]_c=0, \qquad [z_{i},y_{jk}]_c=0, & &i<j, \ i\le k\in\I_{\theta-2};
\\
& [y_{ij},\widetilde{z}_{k}]_c=0, \qquad [\widetilde{z}_{i},y_{jk}]_c=0, & &i<j, \ i\le k\in\I_{\theta-2};
\\
& [y_{ij},y_{kl}]_c=0, & &i,j,k,l \in\I_{\theta-2};
\end{aligned}
\\ \label{eq:rels-ZU-Dtheta-1-thetamenos1-PRV1}
&\begin{aligned}
& x_{\theta}^N=0; &
&z_i^N =\widetilde{z}_i^N =0, \ i\in\I_{\theta-2}; &
&y_{ij}^N =0, \ i\le j\in\I_{\theta-1}.
\end{aligned}
\end{align}

\subsubsection{$\Jb = \I_{2, \theta}$}\label{exa:Dtheta-2-theta}

Here, $Z_U$ is spanned by $\widetilde{w}_{\theta}:= [x_{\theta},w_{\theta-2}]$,
\begin{align*}
w_i&:=x_{i \, i-1 \dots 1}, \ i\in\I, & v_j&:=[[x_{\theta},x_{\theta-2\dots j}]_c, w_{\theta-1}]_c, \ j\in\I_{\theta-2}.
\end{align*}

The algebra $\toba(Z_U)$ is  presented by $w_i$, $i\in\I$, $\widetilde{w}_{\theta}$, $v_j$, $j\in\I_{\theta-2}$, with defining relations
\begin{align} \label{eq:rels-ZU-Dtheta-2-theta-qsr}
&\begin{aligned}
& [w_i,w_j]_c=0, \ i>j\in\I; \quad [v_i,v_j]_c=0, \ i>j\in\I_{\theta-2};
\\
&[\widetilde{w}_{\theta}, w_i]_c = 0, \ i\in\I; \quad [\widetilde{w}_{\theta}, v_i]_c = 0, \ i\in\I_{\theta-2};
\\ 
& [v_j,w_{j-1}]_c = \bq_{\alpha_{(j\theta-2)}+\alpha_{\theta}, \alpha_{(1 \ \theta-1)}}(q-1) w_{\theta-1}\widetilde{w}_{\theta}, \quad j\in\I_{\theta-2}; 
\\ 
&[v_j,w_i]_c=0, \quad i\in\I, j\in\I_{\theta-2}-\{i+1\};
\end{aligned}
\\ \label{eq:rels-ZU-Dtheta-2-theta-PRV1}
&\begin{aligned}
&w_i^N =0, \ i\in\I; &
& v_j^N =0, \ j\in\I_{\theta-2}; & & \widetilde{w}_{\theta}^N=0.
\end{aligned}
\end{align}

\section{Decompositions with a block}\label{sec:decomp-block}

Below we follow the paper \cite{AAH}.

\subsection{A block and a point, weak interaction}\label{subsec:decomp-block-pt-weak}

Let $W$ be a braided vector space of dimension 3 with braiding given in the  basis $(x_i)_{i\in\I_3}$ by
\begin{align}\label{eq:braiding-block-point}
(c(x_i \otimes x_j))_{i,j\in \I_3} &=
\begin{pmatrix}
\epsilon x_1 \otimes x_1&  (\epsilon x_2 + x_1) \otimes x_1& q_{12} x_3  \otimes x_1
\\
\epsilon x_1 \otimes x_2 & (\epsilon x_2 + x_1) \otimes x_2& q_{12} x_3  \otimes x_2
\\
q_{21} x_1 \otimes x_3 &  q_{21}(x_2 + a x_1) \otimes x_3& q_{22} x_3  \otimes x_3
\end{pmatrix}.
\end{align}
Let $V$,  respectively $U$, be the subspace generated by $x_1, x_2$,
respectively $x_3$; $V$ is a block and $U$ is a point. 
The scalar $q_{12}q_{21}$ is called the \emph{interaction} between the block and the point.
As in \cite{AAH}, the \emph{ghost} is 
$\ghost = \begin{cases} -2a, &\epsilon = 1, \\
a, &\epsilon = -1.
\end{cases}$ If $\ghost \in \N$, then we say that the ghost is \emph{discrete}.

\begin{theorem}\label{th:aah-block-point} \cite[4.1]{AAH}
$\toba(W)$ has finite Gelfand-Kirillov dimension if and only if $\epsilon$, 
$q_{22}$, the interaction and the ghost are as in \cite[Table 5]{AAH}.
\end{theorem}

From now on,  we assume that $\GK \toba(W) < \infty$, i.e. that it is as in \cite[Table 5]{AAH}.
Our aim is to compute $Z_V = \ad_c\toba (U) (V)$.
Clearly,  $\GK \toba(Z_V) = \begin{cases} \GK \toba(W) - 1, &q_{22} = 1, \\
\GK \toba(W), &q_{22}  \neq 1.
\end{cases}$

We first deal with the case when the interaction is weak, i.e. that $q_{12}q_{21}= 1$.
Let $y_i = x_i$, $i\in \I_2$, and
\begin{align}\label{eq:zt}
y_{i + 2} &= (\ad_c x_3)^i x_2, & i &\in \N.
\end{align}

\begin{lema}\begin{enumerate}
		\item If $q_{22} = \pm 1$, then $(y_{i})_{i \in \I_3}$ is a basis of  $Z_V$, with braiding
		$(c(y_i \otimes y_j))_{i,j\in \I_3}=$
\begin{align}\label{eq:braiding-block-pt-weak-1}
 &
\begin{pmatrix}
\epsilon y_1 \otimes y_1&  (\epsilon y_2 + y_1) \otimes y_1& \epsilon q_{12} y_3  \otimes y_1
\\
\epsilon y_1 \otimes y_2 & (\epsilon y_2 + y_1) \otimes y_2& \epsilon q_{12} y_3  \otimes y_2
\\
\epsilon q_{21} y_1 \otimes y_3 &  \epsilon q_{21}(y_2 + (a + \epsilon) y_1) \otimes y_3 - a \epsilon y_3 \otimes y_1 & \epsilon q_{22} y_3  \otimes y_3
\end{pmatrix}.
\end{align}

\item If $q_{22} \in \G'_3$, then $(y_{i})_{i \in \I_4}$ is a basis of  $Z_V$, with braiding  given by
\eqref{eq:braiding-block-pt-weak-1} when $i,j\in \I_3$, except that
\begin{align}\label{eq:braiding-block-pt-weak-omega-33}
c(y_3 \otimes y_3) &= \epsilon q_{22} y_3  \otimes y_3 - \epsilon aq_{12} y_4  \otimes y_1;
\end{align}
\begin{align}\label{eq:braiding-block-pt-weak-omega}
 c(y_i \otimes y_4) &=
\begin{pmatrix}
q_{12}^2 y_4 \otimes y_1&  q_{12}^2 y_4 \otimes y_2 & q_{12}q_{22}^2 y_4 \otimes y_3 & q_{22} y_4 \otimes y_4 \end{pmatrix}.
\end{align}
\begin{align}\label{eq:braiding-block-pt-weak-omega-2}
\begin{aligned}
&c(y_4 \otimes y_i) =
\\ &\begin{pmatrix}
q_{21}^2 y_1 \ot y_4 
\\
q_{21}^2 \left(y_2+ (2a+1)y_1\right) \ot y_4 + (1-q_{22}^2) q_{21} y_3 \ot y_3 + a(q_{22}-1) y_4 \ot y_1
\\
q_{21}q_{22}^2 y_3 \ot y_4 + (q_{22}-1) y_4 \ot y_3
\end{pmatrix}.
\end{aligned}\end{align}
\end{enumerate}
\end{lema}

\pf First, $\ad_c(x_3) x_1 =0$ because the interaction is weak. Thus $Z_V$ is generated by $y_i$, $i\in \N$. Observe that $y_i \in T^{i-1}(W)$ when $i \geq2$; 
thus the non-zero $y_i$'s are linearly independent.
Notice that
\begin{align*}
y_3 &= x_3x_2 - q_{21} (x_2 +a x_1)x_3, & &y_4 = x_3^2x_2 - q_{21} (2)_{q_{22}} x_3x_2x_3 + q_{21}^2 q_{22} x_2x_3^2 
\\ &&&- q_{21} a (2)_{q_{22}} x_3x_1x_3 + 2 q_{21}^2 a  q_{22} x_1x_3^2.
\end{align*}
Observe that $\partial_1(y_3) = -ax_3 \neq 0$, hence $y_3 \neq 0$. Also,
\begin{align*}
\partial_1(y_4) &= a(q_{22} - 1) x_3^2,& \partial_2(y_4) &= 0 = \partial_3(y_4).
\end{align*}
Hence, if $q_{22} = \pm 1$, then $y_4=0$, since $x_3^2 = 0$ when $q_{22} =-1$.  If 
$q_{22} \in \G'_3$, then $x_3^2 \neq 0$ but $x_3^3 = 0$, thus $y_4 \neq 0$ and $y_5=0$.
Now 
\begin{align*}
\delta(y_1) &= g_1 \otimes y_1,&\delta(y_2) &= g_1 \otimes y_2,& \delta(y_3) &= g_1g_2 \otimes y_3 -a x_3g_1 \otimes y_1.
\end{align*}
From here \eqref{eq:braiding-block-pt-weak-1} and \eqref{eq:braiding-block-pt-weak-omega-33} follow using that
\begin{align*}
g_1 \cdot y_3 &= \epsilon q_{12} y_3,& g_2 \cdot y_3 &=  q_{22}q_{21} y_3.
\end{align*}
We also compute
\begin{align*}
&\delta(y_4) = g_1g_2^2 \otimes y_4 + (1-q_{22}^2) x_3 g_1g_2 \otimes y_3 +a(q_{22}-1) x_3^2 g_1 \otimes y_1;
\\
&g_1 \cdot y_4= q_{12}^4 y_4, \qquad g_2 \cdot y_4 =  q_{22}^2q_{21} y_4.
\end{align*}
Now \eqref{eq:braiding-block-pt-weak-omega}, \eqref{eq:braiding-block-pt-weak-omega-2} follow by direct computation.
\epf

Recall that the defining relation of the Jordan plane is
\begin{align}
\label{eq:rels-Jordan} &x_2x_1-x_1x_2+\frac{1}{2}x_1^2,
\end{align}
while for the super Jordan plane the defining relations are
\begin{align}
\label{eq:rels-superJordan} &x_1^2, & & x_2x_{21}- x_{21}x_2 - x_1x_{21}.
\end{align}

We introduce the elements
\begin{align}\label{eq:ztilde-defn}
\zt_t&:=(\ad_c y_2)^{t-1}y_3, & &t\in\N.
\end{align}
They are related with the elements $z_t$ for the Nichols algebras in \cite[\S 4]{AAH}:
\begin{align}\label{eq:ztilde-z-relation}
\zt_t&:=-\epsilon^tq_{12} z_t-\delta_{1,t} aq_{12} \, x_1x_3, & \mbox{for all }&t\in\N.
\end{align}

\subsubsection{Case $V=\lstr(1,\ghost)$, $\ghost \in \N$}\label{subsubsection:lstr-11disc}

\begin{prop} \label{pr:lstr-11disc} The algebra 	$\toba(Z_V)$ is presented by generators $y_1,y_2,y_3$ and relations \eqref{eq:rels-Jordan},
\begin{align}
y_1y_3&=q_{12} \, y_3y_1,  \label{eq:lstr-rels&11disc-1} 
\\
(\ad_c y_3)^2 y_2&=0,  \label{eq:lstr-rels&11disc-qserre} 
\\
\zt_{\ghost+1}&=0,  \label{eq:lstr-rels&11disc-qserre-2} 
\\
\zt_t \zt_{t+1}&=q_{12}^{-1} \, \zt_{t+1}\zt_t, & 1\le & t < \ghost. \label{eq:lstr-rels&11disc-2}
\end{align}
$\toba(Z_V)$ has a PBW-basis
\begin{align*}
B=\{ y_1^{m_1} y_2^{m_2} \zt_{\ghost}^{n_{\ghost}} \dots \zt_1^{n_1}: m_i, n_j \in\N_0\};
\end{align*}
hence $\GK \toba(Z_V) = 2+\ghost$.
\end{prop}

\pf 
Note that $y_i=x_i$, $i=1,2$, determine a braided vector subspace of Jordan type, so \eqref{eq:rels-Jordan} holds in $\toba(Z_V)$, while \eqref{eq:lstr-rels&11disc-1} and \eqref{eq:lstr-rels&11disc-qserre} by direct computation.
Relations \eqref{eq:lstr-rels&11disc-qserre-2} and  \eqref{eq:lstr-rels&11disc-2} are 0 in $\toba(Z_V)$ by \eqref{eq:ztilde-z-relation} and \cite[Lemma 4.13]{AAH}.
Hence the quotient $\widetilde{\toba}$ of $T(V)$
by \eqref{eq:rels-Jordan}, \eqref{eq:lstr-rels&11disc-1}, \eqref{eq:lstr-rels&11disc-qserre}, \eqref{eq:lstr-rels&11disc-qserre-2} and  \eqref{eq:lstr-rels&11disc-2} projects onto $\toba(Z_V)$.

We claim that the subspace $I$ spanned by $B$ is a right ideal of $\widetilde{\toba}$. The proof follows as in \cite[Proposition 4.16]{AAH}. 
As $1\in I$, $\widetilde{\toba}$ is spanned by $B$.

To prove that $\widetilde{\toba} \simeq \toba(Z_V)$, it remains to show that
$B$ is linearly independent in $\toba(Z_V)$. This follows from the decomposition \eqref{eq:braided-bosonization-intro}, i.e.
\begin{align*}
\toba(\lstr(1,\ghost)) \simeq \toba(Z_V) \# \Bbbk [x_3]
\end{align*}
and \cite[Proposition 4.16]{AAH}. Then $B$ is a basis of $\toba(Z_V)$ and $\widetilde{\toba}=\toba(Z_V)$.
The computation of $\GK$ follows from the Hilbert series at once.
\epf

\begin{theorem}\label{thm:laistrygonian-ZV}
	The algebra $\toba(Z_V)$ is an iterated Ore extension; thus it is strongly noetherian, AS-regular and Cohen-Macaulay domain.
\end{theorem}

\pf
Analogous to Theorem \ref{thm:laistrygonian}.
\epf

\subsubsection{Case $V=\lstr(-1,\ghost)$, $\ghost \in \N$}\label{subsubsection:lstr-1-1disc}

\begin{prop} \label{pr:lstr1-1disc}
The algebra $\toba(Z_V)$ is presented by generators $y_1,y_2, y_3$ and relations \eqref{eq:rels-Jordan}, 
\eqref{eq:lstr-rels&11disc-1}, \eqref{eq:lstr-rels&11disc-qserre} and 
\begin{align}\label{eq:lstr-rels&1-1disc}
\zt_t^2&=0, & 1\le& t\le \ghost.
\end{align}
The set
\begin{align*}
B=\{ x_1^{m_1} x_2^{m_2} \zt_{\ghost}^{n_{\ghost}} \dots \zt_1^{n_1}: n_i \in\{0,1\}, m_j \in\N_0 \}
\end{align*}
is a basis of $\toba(Z_V)$ and $\GK \toba(Z_V) = 2$.
\end{prop}

\pf
Analogous to Proposition \ref{pr:lstr-11disc}.
\epf

\subsubsection{Case $V=\lstr_-(1,\ghost)$, $\ghost \in \N$}\label{subsubsection:lstr--11disc}

\begin{prop} \label{pr:lstr--11disc} The algebra 	$\toba(Z_V)$ is presented by generators $y_1,y_2, y_3$ and relations \eqref{eq:rels-superJordan}, \eqref{eq:lstr-rels&11disc-1} and
\begin{align}
\zt_{1+2\ghost}&=0, \label{eq:lstr-rels&-11disc-1} \\
y_{21}\zt_1& = q_{12}^2 \, \zt_1 y_{21},  \label{eq:lstr-rels&-11disc-2} \\
\zt_{2k+1}^2&=0, &  0\le & k < \ghost, \label{eq:lstr-rels&-11disc-3} \\
\zt_{2k} \zt_{2k+1}&= q_{12}^{-1} \, \zt_{2k+1}\zt_{2k}, & 1\le & k < \ghost. \label{eq:lstr-rels&-11disc-4}
\end{align}
The set
\begin{align*}
B=\{ y_1^{m_1} y_{21}^{m_2} y_2^{m_3} \zt_{2\ghost}^{n_{2\ghost}} \dots \zt_1^{n_1}: m_1, n_{2k+1} \in\{0,1\}, m_2, m_3, n_{2k} \in\N_0 \}
\end{align*}
is a basis of $\toba(Z_V)$ and $\GK \toba(Z_V) = \ghost+2$.
\end{prop}

\pf
Analogous to Proposition \ref{pr:lstr-11disc}.
\epf

\subsubsection{Case $V=\lstr_-(-1,\ghost)$, $\ghost \in \N$}\label{subsubsection:lstr--1-1disc}

\begin{prop} \label{pr:lstr-1-1disc} The algebra $\toba(Z_V)$ is presented by generators $y_1,y_2, y_3$ and relations \eqref{eq:rels-superJordan}, \eqref{eq:lstr-rels&11disc-1}, \eqref{eq:lstr-rels&11disc-qserre}, \eqref{eq:lstr-rels&-11disc-1}, \eqref{eq:lstr-rels&-11disc-2} and
\begin{align}
\zt_{2k}^2&=0, &  1\le & k \le \ghost, \label{eq:lstr-rels&-1-1disc-1} \\
\zt_{2k-1} \zt_{2k}&= -q_{12}^{-1} \zt_{2k}\zt_{2k-1}, & 0< & k \le \ghost. \label{eq:lstr-rels&-1-1disc-2}
\end{align}
The set
\begin{align*}
B=\{ y_1^{m_1} y_{21}^{m_2} y_2^{m_3} \zt_{2\ghost}^{n_{2\ghost}} \dots \zt_1^{n_1}: m_1, n_{2k} \in\{0,1\}, m_2, m_3, n_{2k-1} \in\N_0 \}
\end{align*}
is a basis of $\toba(Z_V)$ and $\GK \toba(Z_V) = \ghost+2$.
\end{prop}

\pf
Analogous to Proposition \ref{pr:lstr-11disc}.
\epf

\subsubsection{Case $V=\lstr(\omega,1)$, $\omega \in \G'_3$} \label{subsubsection:lstr-omega-1}

\begin{prop} \label{pr:lstr1omega1} The algebra $\toba(Z_V)$ is presented by generators $y_1,y_2,y_3, y_4$ and relations  
\eqref{eq:rels-Jordan},
\eqref{eq:lstr-rels&11disc-1} and
\begin{align}
\label{eq:rels-lstr1omega1-1}
y_1 y_4 &= q_{12}^2 y_4 y_1,
\\
\label{eq:rels-lstr1omega1-2}
y_2 y_3 &= q_{12} y_3 y_2,
\\
\label{eq:rels-lstr1omega1-3}
y_4 y_2 &= q_{21}^2 y_2 y_4 + q_{21}(1-\omega) \, y_3^2,
\\
\label{eq:rels-lstr1omega1-4}
y_3 y_4 &= q_{12}\omega^2 y_4 y_3,
\\
\label{eq:rels-lstr1omega1-5}
y_3^3 &= y_4^3=0.
\end{align}
The set
\begin{align*}
B=\{ y_1^{m_1} y_2^{m_2} y_3^{n_1} y_4^{n_2}: m_i\in\N_0, 0 \le n_j\le 2 \}
\end{align*}
is a basis of $\toba(Z_V)$ and $\GK \toba(Z_V) = 2$.
\end{prop}

\pf
Relations \eqref{eq:rels-Jordan} and \eqref{eq:lstr-rels&11disc-1} are 0 in $\toba(Z_V)$ as in Proposition \ref{pr:lstr-11disc}. Now \eqref{eq:rels-lstr1omega1-1}-\eqref{eq:rels-lstr1omega1-5} follow from \cite[Lemmas 4.23 \& 4.24]{AAH}. Hence the quotient $\widetilde{\toba}$ of $T(V)$
by all these relations projects onto $\toba(Z_V)$. Since the subspace $I$ spanned by $B$ is a right
ideal of $\widetilde{\toba}$ and $1\in I$, $\widetilde{\toba}$ is spanned by $B$.
To prove that $\widetilde{\toba} \simeq \toba(Z_V)$, it remains to show that
$B$ is linearly independent in $\toba(Z_V)$. 
This follows from the decomposition 
\begin{align*}
\toba(\lstr(\omega,1)) \simeq \toba(Z_V) \# \Bbbk [x_3]
\end{align*}
as in \eqref{eq:braided-bosonization-intro}
and \cite[Proposition 4.25]{AAH}. Then $B$ is a basis of $\toba(Z_V)$, $\widetilde{\toba}=\toba(Z_V)$ and $\GK \toba(Z_V)=2$.
\epf

\subsection{A block and a point, mild interaction}\label{subsec:decomp-block-pt-mild}
Here we keep the notation as in the previous Subsection but we assume that the interaction is
mild, that is $q_{12}q_{21}= -1$. We consider the unique Nichols algebra of finite $\GK$, called the Cyclop Nichols algebra: here, $\epsilon = q_{22}=-1$. Let 
\begin{align}\label{eq:mild-yt-def}
y_i &= x_i, & y_{i + 2} &= (\ad_c x_3) x_i, & i & \in \I_2.
\end{align}

\begin{lema}
A basis of  $Z_V$ is given by $(y_{i})_{i \in \I_4}$, with braiding

\begin{align}\label{eq:braiding-block-pt-mild-1}
c(y_i\otimes y_j) & = 
\begin{cases}
-y_1\otimes y_i, & j=1, \\
(y_1-y_2)\otimes y_i, & j=2, \\
-q_{12}y_3\otimes y_i, & j=3, \\
q_{12}(y_3-y_4)\otimes y_i, & j=4,
\end{cases}, \qquad i=1,2;
\\
\label{eq:braiding-block-pt-mild-2}
c(y_3\otimes y_j) & = 
\begin{cases}
-q_{21}y_1\otimes y_3-2y_3\otimes y_1, & j=1, \\
-q_{21}y_2\otimes y_3+2(y_3-y_4)\otimes y_1, & j=2, \\
-y_3\otimes y_3, & j=3, \\
-y_4\otimes y_3, & j=4;
\end{cases}
\\
\label{eq:braiding-block-pt-mild-3}
c(y_4 \otimes y_j) & = 
\begin{cases}
-q_{21}y_1\otimes y_4-2y_3\otimes y_2-y_3\otimes y_1, & j=1, \\
-q_{21}y_2\otimes y_4+2(y_3-y_4)\otimes (2y_2+y_1), & j=2, \\
-y_3\otimes y_4, & j=3, \\
-y_4\otimes y_4, & j=4.
\end{cases}
\end{align}
\end{lema}

\pf First, $(\ad_c x_3)^2 x_i =0$ because $x_3^2=0$. Thus $Z_V$ is generated by $y_i$, $i\in \I_4$. We claim that the $y_i$'s are linearly independent. Indeed,
\begin{align*}
\partial_1(y_3) &= 2x_3,
&
\partial_2(y_3) &= 0,
& 
\partial_1(y_4) &= x_3,
&
\partial_2(y_4) &= 2x_3,
\end{align*}
and $\partial_i(y_j)=\delta_{ij}$ for $i,j\in\I_2$.
Now 
\begin{align*}
\delta(y_1) &= g_1 \otimes y_1, & 
\delta(y_3) &= g_1g_2 \otimes y_3 +2 x_3g_1 \otimes y_1,
\\
\delta(y_2) &= g_1 \otimes y_2, & 
\delta(y_4) &= g_1g_2 \otimes y_3 +2 x_3g_1 \otimes y_2  + x_3g_1 \otimes y_1.
\end{align*}
From here \eqref{eq:braiding-block-pt-mild-1}, \eqref{eq:braiding-block-pt-mild-2} and \eqref{eq:braiding-block-pt-mild-3}  follow by direct computation.
\epf

We set $y_{14}=(\ad_c y_1)y_4=y_1y_4+q_{12}y_4y_1-q_{12}y_3y_1$.

\begin{prop} \label{pr:cyclop} The algebra 	$\toba(Z_V)$ is presented by generators $(y_j)_{j\in \I_4}$ and relations \eqref{eq:rels-superJordan},
\begin{align}\label{eq:block-cyclop-rels-1}
y_1y_3 + q_{12}y_3y_1 &= 0, & y_2y_3+q_{12}y_3y_2 &= -q_{12} y_{14} - q_{12}y_3y_1, 
\\ \label{eq:block-cyclop-rels-2}
y_2y_4 + q_{12} y_4y_2 &=0, & y_3y_4+y_4y_3 &= 0,  
\\ \label{eq:block-cyclop-rels-3}
y_{14}^2 &=0, & y_3^2 &=0, \qquad y_4^2=0.
\end{align}
The set
\begin{align*}
B=\{ y_1^{m_1} y_{21}^{m_2} y_2^{m_3} y_{14}^{n_1} y_3^{n_2} y_4^{n_3}: m_1, n_{i} \in\{0,1\}, m_2, m_3, \in\N_0 \}
\end{align*}
is a basis of $\toba(Z_V)$ and $\GK \toba(Z_V) = 2$.
\end{prop}

\pf
All the quadratic relations belong to $\ker (\id + c)$, the quantum symmetrizer of degree 2, hence they are defining relations of $Z_V$. Now the first relation of \eqref{eq:block-cyclop-rels-3} follows from \cite[Lemma 4.34]{AAH}.
Hence the quotient $\widetilde{\toba}$ of $T(V)$
by \eqref{eq:block-cyclop-rels-1}, \eqref{eq:block-cyclop-rels-2} and \eqref{eq:block-cyclop-rels-3} projects onto $\toba(Z_V)$. Using these relations we check that the subspace $I$ spanned by $B$ is a right
ideal of $\widetilde{\toba}$. Since $1\in I$, $\widetilde{\toba}$ is spanned by $B$.

To prove that $\widetilde{\toba} \simeq \toba(Z_V)$, it remains to show that
$B$ is linearly independent in $\toba(Z_V)$. 
This follows from the decomposition 
\begin{align*}
\toba(\lstr(\omega,1)) \simeq \toba(Z_V) \# \Bbbk [x_3]
\end{align*}
as in \eqref{eq:braided-bosonization-intro}
and \cite[Proposition 4.39]{AAH}. Then $B$ is a basis of $\toba(Z_V)$, $\widetilde{\toba}=\toba(Z_V)$ and $\GK \toba(Z_V)=2$.
\epf

\smallbreak\subsubsection*{Acknowledgements} 
This paper grew from conversations following a talk by Oscar M\'arquez on joint work in progress 
with Dirceu Bagio and Gast\'on A. Garc\'\i a
at the Colloquium Quantum 17 hosted by the University of Talca   (Chile). 
We thank them for sharing their results as well as Mar\'\i a Ronco and Mar\'\i a In\'es Icaza for hospitality.
We also thank Hiroyuki Yamane for pointing out to us the reference \cite{U}.

We are grateful to C.D. Ward and H. West (University of Miskatonic, Arkham)
for pointing out to us a mistake in the proof of Lemma 3.6.

The main results of this paper were communicated at the \emph{XXII Coloquio Latinoamericano de \'Algebra} (Quito, August 2017); the \emph{Reuni\'on Anual de la Uni\'on Matem\'atica Argentina} (Buenos Aires, December 2017);
the Workshop \emph{M\'etodos Categ\'oricos en \'Algebras de Hopf} (Maldonado, December 2017);
the Workshop \emph{Tensor categories, Hopf algebras and quantum groups} (Marburg, January 2018).


\begin{thebibliography}{AAH}
\bibitem[A1]{A-chicago} 
Andruskiewitsch, N.: \emph{Some remarks on Nichols algebras}. In "Hopf algebras", Bergen, Catoiu and Chin (eds.), 25--45, M. Dekker, (2004).

\bibitem[A2]{icm} Andruskiewitsch, N.:
On finite-dimensional Hopf algebras. Proceedings of the ICM Seoul 2014 Vol. II,  117--141 (2014)

\bibitem[A3]{A-leyva} Andruskiewitsch, N.: \emph{An Introduction to Nichols Algebras}. In Quantization, Geometry and Noncommutative Structures in Mathematics and Physics. 
A. Cardona, P. Morales, H. Ocampo, S. Paycha, A. Reyes, eds., pp. 135--195, Springer (2017).

\bibitem[AA1]{AA} Andruskiewitsch, N., Angiono, I.: \emph{On Nichols algebras with generic braiding}. In  Modules and Comodules, Trends in Mathematics. Brzezinski, T.; Gomez
Pardo, J.L.; Shestakov, I.; Smith, P.F. (Eds.), pp. 47--64 (2008). ISBN: 978-3-7643-8741-9.


\bibitem[AA2]{AA-diag-survey} Andruskiewitsch, N., Angiono, I.: \emph{On Finite dimensional Nichols algebras of diagonal type}. 
Bull. Math. Sci. \textbf{7} 353--573 (2017). 



\bibitem[AAH1]{AAH} Andruskiewitsch, N.,  Angiono, I., Heckenberger, I.:
\emph{On finite GK-dimensional Nichols algebras over abelian groups}. Mem. Amer. Math. Soc.,
to appear.



\bibitem[AAH2]{AAH-infinite} Andruskiewitsch, N.,  Angiono, I., Heckenberger, I.:
\emph{On Nichols algebras of infinite rank with finite Gelfand-Kirillov dimension}. \texttt{arXiv:1805.12000}.




\bibitem[AC]{AC} Andruskiewitsch, N., Cuadra, J.: \emph{On the structure of (co-Frobenius) Hopf algebras}. J.
Noncommut. Geom. \textbf{7}  83--104 (2013).






\bibitem[AHS]{AHS} Andruskiewitsch, N., Heckenberger, I., Schneider, H.-J.: 
\emph{The Nichols algebra of a semisimple Yetter-Drinfeld module}, Amer. J. Math. \textbf{132} 1493--1547 (2010).



\bibitem[AS1]{AS-98} Andruskiewitsch, N.,   Schneider, H.-J.,  Lifting of quantum linear spaces and pointed Hopf algebras of order $ p^3$, 
\emph{J. Algebra} \textbf{209} (1998), 658--691.

\bibitem[AS2]{AS-adv} Andruskiewitsch, N.,  Schneider, H.-J.: \emph{Finite quantum groups and Cartan matrices}, Adv. Math. \textbf{154} 1--45 (2000).

\bibitem[AS3]{AS-cambr} Andruskiewitsch, N.,  Schneider, H.-J.: 
\emph{Pointed Hopf algebras}. In Recent developments in Hopf algebras Theory, 
MSRI Publ. \textbf{43} 1--68, Cambridge Univ. Pr. (2002).

\bibitem[AS3]{AS-ann}  Andruskiewitsch, N., Schneider, H.-J.:
On the classification of finite-dimensional pointed Hopf
algebras, Ann. Math. \textbf{171},  375--417  (2010)

\bibitem[An1]{A-jems} Angiono, I.: \emph{A presentation by generators and relations of Nichols algebras of diagonal type and convex orders on root systems}. 
J. Eur. Math. Soc.  \textbf{17} 2643--2671 (2015).


\bibitem[An2]{A-presentation} Angiono, I.: 
\emph{On Nichols algebras of diagonal type}. 
J. Reine Angew. Math.  \textbf{683}  189--251 (2013).

\bibitem[An3]{A-pre-Nichols}  Angiono, I.: 
\emph{Distinguished pre-Nichols algebras}. 
Transf. Groups \textbf{21} 1--33 (2016).

\bibitem[AnG]{AG} Angiono, I.,  Garc\'\i a Iglesias, A.:
\emph{Liftings of Nichols algebras of diagonal type II. All liftings are cocycle deformations}. \texttt{arXiv:1605.03113}.

\bibitem[AnG2]{AG-survey} Angiono, I.,  Garc\'\i a Iglesias, A.:
\emph{Pointed Hopf algebras: a guided tour to the liftings}. \texttt{arXiv:1807.07154}.

\bibitem[ArS]{ArS} Artin,  M., Schelter,  W. F.: 
\emph{Graded algebras of global dimension 3}. Adv. Math.
\textbf{66} 171--216 (1987).

\bibitem[AST]{AST} Artin,  M., Schelter,  W. F., Tate, J.: 
\emph{Quantum deformations of $GL_n$}.
Comm. Pure Appl. Math. \textbf{44} 879--895 (1991).

\bibitem[ASZ]{ASZ} Artin,  M., Small, L. W. , Zhang,  J. J.: 
\emph{Generic flatness for strongly Noetherian algebras}. J. Algebra \textbf{221}  579--610 (1999).


\bibitem[CL]{CL}
Cuntz, M., Lentner, S.: \emph{A simplicial complex of Nichols algebras}. Math. Z. \textbf{285} 647--683 (2017).

\bibitem[DoT]{DT} Doi, Y., Takeuchi, M.:
\emph{Multiplication alteration by two-cocycles - the quantum version}.
Commun. Alg. {\bf 22} 5715--5732 (1994).

\bibitem[Dr1]{Dr} Drinfeld, V. G.:
\emph{Hopf algebras and the quantum Yang-Baxter equation}.
Sov. Math. Dokl. \textbf{32} 256--258 (1985).

\bibitem[E]{E} Elle, S.: \emph{Classification of relation types of Ore extensions of dimension 5}. 
Commun. Alg. \textbf{45}  1323--1346 (2017).


\bibitem[GG]{GGi} Garc\'{\i}a, G. A.,  Giraldi, J. M. J.: 
\emph{On Hopf Algebras over quantum subgroups}. 
 J. Pure Appl. Algebra \textbf{223} (2019), 738--768.

\bibitem[GM]{GM} 
Garc\'{\i}a, G. A., Mastnak, M.: 
\emph{Deformation by cocycles of pointed Hopf algebras over non-abelian groups}. 
Math. Res. Lett. \textbf{22} 59--92 (2015). 


\bibitem[Gr]{gr-jalg} Gra\~na, M.: 
\emph{A Freeness Theorem for Nichols Algebras}. 
J. Algebra \textbf{231} 235--257 (2000).

\bibitem[H1]{H-inv} Heckenberger, I.: 
\emph{The Weyl groupoid of a Nichols algebra of diagonal type}. Invent. Math. \textbf{164} 175--188 (2006).

\bibitem[H2]{H-classif RS} Heckenberger, I.: \emph{Classification of arithmetic root
	systems}. Adv. Math. \textbf{220}  59--124 (2009).

\bibitem[HS1]{HS-coideal subalg} Heckenberger, I., Schneider, H.-J.: 
\emph{Right coideal subalgebras of Nichols algebras and the Duflo order on the Weyl groupoid}.  
Israel J. Math. \textbf{197} 139--187 (2013).

\bibitem[HS2]{HS-adv} Heckenberger, I., Schneider, H.-J.:
\emph{Yetter-Drinfeld modules over bosonizations of dually paired Hopf algebras}. 
 Adv. Math. \textbf{244} 54--394 (2013).


\bibitem[HV]{HV-rank>2} Heckenberger, I.,  Vendramin,  L.:
\emph{A classification of Nichols algebras of semi-simple Yetter-Drinfeld modules over non-abelian groups}. 
J. Eur. Math. Soc. \textbf{19} 299--356 (2017). 

\bibitem[HY1]{HY-Shap-det} Heckenberger, I., Yamane, H.: 
\emph{Drinfel'd doubles and Shapovalov determinants}. 
Rev. Un. Mat. Argentina \textbf{51} 107--146 (2010).

\bibitem[HY2]{HY-Weyl-gpd} Heckenberger, I., Yamane, H.: 
\emph{A generalization of Coxeter groups, root systems, and Matsumoto's theorem}. 
Math. Z.
\textbf{259} 255--276 (2008).

\bibitem[HX]{HX} Hu, N., Xiong,  R.: 
\emph{On families of Hopf algebras without the dual Chevalley property}, Rev. Un. Mat. Argentina,
59 (2018), 443--469.



\bibitem[L]{Lu} Lusztig, G.:\emph{ Introduction to quantum groups}. Birkh\"auser (1993).


\bibitem[LW]{LW} Li, J., Wang, X.:
\emph{Some five-dimensional Artin-Schelter regular algebras obtained by deforming a Lie algebra}.
J. Alg. Appl. \textbf{15} (04) 1650060 (2016).

\bibitem[M]{Majid}
Majid,  S.: \emph{Doubles of quasitriangular Hopf algebras}. Comm. Algebra \textbf{19} 3061--3073 (1991).


\bibitem[Mo]{Mo-libro} Montgomery, S.: 
\emph{Hopf algebras and their actions on rings}, CMBS \textbf{82},  
Amer. Math. Soc. (1993).

\bibitem[NS]{NS} Nevins, T. A.,  Stafford, J. T.: 
\emph{Sklyanin algebras and Hilbert schemes of points}. Adv. Math.,
\textbf{210} 405--478 (2007).


\bibitem[PV]{PV} Pogorelsky, B., Vay, C.: 
\emph{Verma and simple modules for quantum groups at non-abelian groups}. Adv. Math. \textbf{301} 423--457 (2016).

\bibitem[RS]{RS} Radford, D. E., Schneider, H.-J.: 
\emph{On the simple representations of generalized
quantum groups and quantum doubles}. J. Algebra \textbf{319} 3689--3731 (2008).


\bibitem[R]{R}  Rosso, M.: \emph{Quantum groups and quantum shuffles}. Invent. Math. \textbf{133} 399--416 (1998).

\bibitem[S]{S} Schauenburg, P.: \emph{Hopf bi-Galois extensions}. Comm. Algebra \textbf{24} 3797--3825 (1996).



\bibitem[U1]{U1}  Ufer, S.: 
\emph{PBW bases for a class of braided Hopf algebras}. 
J. Algebra \textbf{280} 84--119 (2004). 

\bibitem[U2]{U2}  Ufer, S.: 
\emph{Triangular braidings and pointed Hopf algebras}.
J. Pure Appl. Algebra \textbf{210} 307--320 (2007). 


\bibitem[U3]{U}  Ufer, S.: 
\emph{Braided Hopf algebras of triangular type}. PhD thesis (2004). \newline
\texttt{https://edoc.ub.uni-muenchen.de/2477/1/ufer\char`_stefan.pdf}

\bibitem[X1]{X1} Xiong,  R.: \emph{On Hopf algebras over the unique 12-dimensional Hopf algebra without the dual Chevalley property}. Commun. Algebra, to appear.

\bibitem[X2]{X2} Xiong,  R.: \emph{Finite-dimensional Hopf algebras over the smallest non-pointed basic Hopf algebra}. \texttt{arXiv:1801.06205}.


\bibitem[X3]{X3} Xiong,  R.: \emph{On Hopf algebras over basic Hopf algebras of dimension 24}. \texttt{arXiv:1809.03938}.

\bibitem[WW]{WW} Wang Q., Wu Q. S.: \emph{A class of AS-regular algebras of dimension five}, J. Algebra \textbf{362} (2012), 117--144.

\bibitem[ZZ]{ZZ} Zhang,  J. J., Zhang, J.: 
\emph{Double extension regular algebras of type}. 
J. Algebra \textbf{322} 373--409 (2009).

\end{thebibliography}
\end{document}